\documentclass[11pt,twoside]{amsart}

\usepackage{amssymb}
\usepackage{amsxtra}
\usepackage{graphicx}
\usepackage{dcpic,pictexwd}

\theoremstyle{plain}
\newtheorem{thm}{Theorem}
\newtheorem{theorem}{Theorem}[section]
\newtheorem{prop}[theorem]{Proposition}

\newtheorem{defn}[theorem]{Definition}

\newtheorem{lemma}{Lemma}
\newtheorem{cor}{Corollary}

\newcommand{\al}{\alpha}
\newcommand{\be}{\beta}
\newcommand{\ga}{\gamma}
\newcommand{\de}{\delta}

\newcommand{\Q}{\mathbb{Q}}
\newcommand{\F}{\mathbb{F}}
\newcommand{\R}{\mathbb{R}}
\newcommand{\T}{\mathbb{T}}
\newcommand{\C}{\mathbb{C}}
\newcommand{\Z}{\mathbb{Z}}

\newcommand{\LL}{\mathbb{L}}

\newcommand\ra{\rightarrow}
\newcommand\lra{\longrightarrow}

\newcommand\seta{\{\alpha_i\}_{i=1}^g}
\newcommand\setb{\{\beta_i\}_{i=1}^g}
\newcommand\setg{\{\gamma_i\}_{i=1}^g}

\newcommand{\Si}{\Sigma}
\newcommand{\Siab}{(\Sigma,\, \seta ,\, \setb; w)}

\newcommand{\Siabg}{(\Sigma,\, \seta,\, \setb,\, \setg; w)}

\newcommand{\hatc}[1]{\widehat{CF}(#1)}

\newcommand{\ben}{\begin{enumerate}}
\newcommand{\een}{\end{enumerate}}

\newcommand{\POz}{P.\ Ozsv{\'a}th\,}
\newcommand{\ZSz}{Z.\ Szab{\'o}\,}

\newcommand{\Sp}[1]{\mathfrak{#1}}

\newcommand{\Sps}[2]{\Sp{#1}_{#2}}

\newcommand{\codeS}[2]{I^{#1} < \cdots < I^{#2}} 
\newcommand{\hth}[2]{\widehat{\Theta}_{#1,\,#2}}

\newcommand{\mir}[1]{\overline{#1}}
\newcommand{\der}{\partial}

\newcommand{\HomT}[1]{H_{\ast}(#1)}
\newcommand{\csquare}[8]{\begindc{\commdiag}[6]
	\obj(0,0)[I1]{$#3$}
	\obj(15,0)[I2]{$#4$}
	\obj(0,10)[I3]{$#1$}
	\obj(15,10)[I4]{$#2$}
	\mor{I1}{I2}{$#8$}
	\mor{I3}{I4}{$#5$}
	\mor{I4}{I2}{$#7$}
	\mor{I3}{I1}{$#6$} \enddc }

\begin{document}
\title{Notes on the Heegaard-Floer Link Surgery Spectral Sequence}

\author[L. Roberts]{Lawrence Roberts}
\address{Department of Mathematics, Michigan State University,
East Lansing, MI 48824}
\email{lawrence@math.msu.edu}
\thanks {The author was supported in part by NSF grant DMS-0353717 (RTG)}
\maketitle

\noindent In \cite{Doub}, \POz and \ZSz constructed a spectral sequence computing the Heegaard-Floer homology $\widehat{HF}(Y_{\LL})$ where $Y_{\LL}$ 
is the result of surgery on a framed link, $\LL$, in $Y$. The terms in the $E^{1}$-page of this spectral sequence are Heegaard-Floer homologies of surgeries on $\LL$ for other framings derived from the original. They used this result to analyze the branched double cover of a link $L \subset S^{3}$ where it was possible to give a simple description of all the groups arising in the $E^{1}$-page. The result is a spectral sequence, over $\F_{2}$, with $E^{2}$ page given by the reduced Khovanov homology of $L$ and converging in finitely many steps to $\widehat{HF}(-\Si(L))$, where $\Si(L)$ is the branched double cover of $S^{3}$ over $L$. Several years later, in \cite{Robe}, \cite{Rob2} adjusted this argument to a setting where the spectral sequence started at refinement of Khovanov homology, and converged to a knot Floer homology. This facilitated the analysis of the knot Floer homology of certain fibered knots. Recently, Olga Plamanevskaya first in \cite{Pla3} used this approach to show that the contact invariant of certain open books was non-vanishing. By generalizing from the double branched cover picture, she and John Baldwin, \cite{Bald}, were able to extend the argument to more general open books. In addition, Eli Grigsby and Stefan Wehrli, \cite{EliS} found a different direction in which to generalize link surgery spectral sequence: to sutured Floer homology. They then showed that the Khovanov categorification of the colored Jones polynomial detects the unknot.\\
\ \\
\noindent This paper aims to prove several foundational results for these types of spectral sequences. First, it reviews in detail the construction of \cite{Doub}, laying out the construction in a manner which will allow the later results in this paper. This occurs in part I, which contains no new results. Part II starts by explaining how to modify part I for knot Floer homology, elaborating on the terse proof given in \cite{Robe}. It then proceeds to the novel part of the paper: explaining how the spectral sequence is invariant of the many choices in its construction, under a suitable equivalence. This forms a prelude to incorporating cobordism maps into the picture. The author included part I, even though there are no new results, because these sections depend heavily on the methods employed in the original construction of the spectral sequence. In the end we will, for example, be able to answer such question as: if $W: Y_{\LL} \ra Y_{\LL,K}$ is a $4$-dimensional oriented cobordism formed by adding a two handle along the framed knot $K \subset Y_{\LL}$, can we find a morphism of the spectral sequences arising from $\LL$ for each end which reflects the cobordism map for $W$? Combined with some of the techniques for analyzing knot Floer homology, one can obtain information about the cobordism maps distinguished by certain $Spin^{c}$-properties. This should form the basis for answering the question about cobordism maps at the end of the introduction to \cite{Doub}. \\
\ \\
\noindent {\bf Convention:} Throughout $Y$ will be a closed, oriented, smooth three manifold. All calculations are assumed to be performed over the finite field $\F_{2}$. No effort has been made to make the signs correct for other rings. Heegaard decompositions use the convention that $(\Si, H_{\al}, H_{\be})$ has $\partial H_{\al} = \Si = - \partial H_{\be}$ with the outward normal first convention. We will generally consider knots and links to lie in $H_{\be}$. Thus, the boundary orientation for the boundary of a knot complement will coincide with the Heegaard surface orientation. Finally, we will always assume, unless otherwise stated, that we are not distiguishing $Spin^{c}$ structures in our Heegaard-Floer homologies, i.e. we always take a direct sum of chain groups across all relevant $Spin^{c}$ structures. \\

\part*{I. The construction of the spectral sequence in \cite{Doub}}
\noindent This part, sections 1-5, recount the original proof of the link surgery spectral given by \POz and \ZSz in \cite{Doub}. We repeat it here in slightly expanded form in order to remind the reader of the notation and conception, and to make more precise some of the choices necessary to the
construction. These sections will be assumed in the second part of the paper, which recounts new results.\\

\section{Triads of Framings}

\begin{center}
\begin{figure}
\includegraphics[scale=0.5]{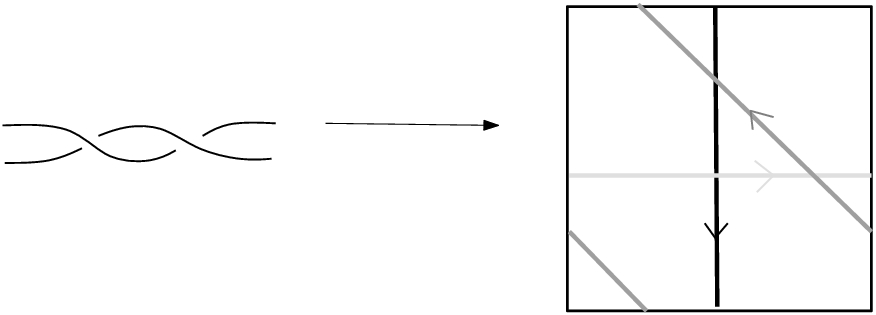}
\caption{A $+1$-framing on the unknot corresponds to the Heegaard triple on the right, where the torus is oriented as a boundary of a tubular neighborhood of the knot. The vertical, black, curve is the meridian and a $\be$ curve; the horizontal curve is the $0$-framing and an $\ga$-curve, and the slanted curve is the $+1$ framing and the $\de$-curve. These are oriented to have $+1$ intersection in the torus. Since this diagram has the opposite orientation of the splitting surface, using the convention in the text, this forms a triad in the Heegaard diagram. Note that there is an $\be\ga\de$-triangle with multiplicity $-1$ in this diagram; when we take a reflection of this diagram, it is a holomorphic triangle.} \label{fig:triad}
\end{figure}
\end{center}

\begin{defn} Let $K$ be a knot in $Y$. A triple of framings for $K$ form a {\em triad} if in $M = Y-int\,N(K)$ there are representatives for the framings thought of as curves, $F_{a}, F_{b}$ and $F_{c}$,  in $\partial M$ such that 
$$
F_{a} \cap F_{b} = F_{b} \cap F_{c} = F_{c} \cap F_{a} = - 1
$$
for algebraic intersections in $\partial M$ oriented as the boundary of $M$. 
\end{defn}
\ \\
{\bf Example:} The $\infty, n, n+1$ framings for a knot in $S^{3}$ form a triad. This is the image of $m$, $l$, and $l + m$ under the map taking $m \ra m$ and $l \ra l + n\cdot m$, so it suffices to verify the property for $\infty, 0, +1$, see Figure \ref{fig:triad}. Note that the outward normal first convention means that the intersection numbers for the standard orientation of the plane are the negatives of those in the boundary. 

\begin{lemma}
Any two oriented framings for a triad determines the third. Furthermore, the roles of the three curves are symmetric under cyclic permutation.
\end{lemma}

\noindent {\bf Proof:} If $a$ is the first oriented curve, $b$ the second oriented curve, and $a \cap b = -1$, then the third curve is $-a-b$. This follows since $a, b$ span $H_{1}(T^{2};\Z)$, and the intersection numbers determine the third oriented curve. Note that had we chosen $a'=b$ and $b'=-a-b$, then $-a'-b'=a$, and had we chosen $a'=-a-b$ and $b'=a$, then $-a'-b'= b$. Furthermore, choosing the opposite orientations, $-a$ and $-b$ on the framings gives $a+b$ for the third framing. As oriented curves, these are different, but as unoriented framings, which is all that matters to Dehn surgery, these are the same. $\Diamond$

\begin{lemma}
There are cobordisms: $W_{1}: Y_{a} \ra Y_{b}$, $W_{2}: Y_{b} \ra Y_{c}$, and $W_{3}: Y_{c} \ra Y_{a}$ each of which consists of a single two handle addition. The compositions $W_{2} \circ W_{1}$, $W_{3} \circ W_{2}$ and $W_{1} \circ W_{3}$ each contain an embedded sphere of self-intersection $-1$.
\end{lemma}

\noindent {\bf Proof:} We describe how to find triads for more general framings than in the example above. We can think of $Y$ as given by a surgery diagram, i.e. an integer framed link which will yield $Y$ when we attach $4$-dimensional two handles according to the given data. Then $K$ can be depicted by a component of a link in $S^{3}$ where the other components are the framed link describing $Y$. This perspective uses, as in the example above, the wrong orientation on $T^{2}$. Thus we will look for three oriented curves with $+1$ intersections in their cycle. Relative to the Seifert framing, the first curve in the triad, $F_{a}$, can be described by a rational number.  To describe the other curves in the triad, we need to convert this framing. That is to say, we need to find a means to think of the result of rational Dehn surgery on $K \subset Y$ as a sequence of integer surgeries (for ease, on $S^{3}$ in addition to the framed link). We do this by attaching a sequence of unknots each linked with its predecessor and the first linked to $K$. First, define
$$
[q_1, q_2, \ldots, q_{m} ]= q_{1} - \frac{1}{q_{2} - \frac{1}{q_{3} - \frac{1}{\ddots - \frac{1}{q_{m}}}}}
$$
Then if we frame $K$ with $q_{1}$ and each link in the chain with the $q_{i}$ in order, the boundary is the equivalent of framing $K$ with 
$[q_{1}, \ldots, q_{m}]$. We rewrite this using the $+$ signs to combine the fractions, which we denote by $(q_{1}, \ldots, q_{m})$. We can convert between these using the formula $[q_{1}, \ldots, q_{m}] = ( q_{1} , - q_{2}, q_{3}, -q_{4}, \ldots, \pm q_{m})$. We denote by $\frac{A_{i}}{B_{i}}$ the fraction generated by only the first $i$ entries $(q_{1}, -q_{2}, q_{3}, \ldots, \pm q_{i})$. There are two well known facts about these fractions: 1) $A_{i+1} = q_{i+1}A_{i} + A_{i-1}$ and $B_{i+1} = q_{i+1}B_{i} + B_{i-1}$ and 2) $A_{i-1}B_{i} - B_{i-1}A_{i} = (-1)^{i-1}$. This last identity represents the algebraic intersection number of the framing curves in the torus boundary, found from the $i-1^{st}$ and $i^{th}$ partial continued fractions. This allows us to append a $p$ framed unknot to the end of the chain setting $q_{m+1} = p$. There are two cases to consider based on the parity of $m$. \\
\ \\
\noindent If $m$ is even, we would have $A_{m+1} = p\,A_{m} + A_{m-1}$ and $B_{m+1} = p\,B_{m} + B_{m-1}$ If we append a second $-1$ framed unknot, we will get $A_{m+2} =  A_{m+1} + A_{m} = (p + 1)\,A_{m} + A_{m-1}$ and $B_{m+2} =(p+1)\,B_{m} +  B_{m-1}$ since the end will be $\ldots, q_{m-1}, -q_{m}, p, +1)$. If we append another  $-1$ framed unknot, we will have $\ldots, -q_{m}, p, +1, -1)$ and $A_{m+3} = - A_{m+2} + A_{m+1} = - A_{m}$ and $B_{m+3} = - B_{m}$. One more produces $A_{m+4} = A_{m+3} + A_{m+2} = p\,A_{m} + A_{m-1} = A_{m+1}$ and $B_{m+4} = B_{m+1}$. Thus we will have the sequence $\frac{A_{m}}{B_{m}},$ $\frac{A_{m+1}}{B_{m+1}},$ $\frac{A_{m+2}}{B_{m+2}}$, $\frac{-A_{m}}{-B_{m}}$. Since $m$ is even, the intersection numbers for these pairs of framings are $+1,-1,+1$. By switching the signs on $(A_{m+2}, B_{m+2})$ we can change the orientation of the framing curve, but not the three manifold found by filling along that curve. This gives the sign pattern $+1,+1,+1$ for the succesive pairs between $m, m+1, m+2$ since the last fraction represents $-1$ times the original framing curve. \\ 
\ \\
\noindent There is, of course, a second possibility where $m$ is odd. Then the end is $\ldots, -q_{m-1}, q_{m}$ $, -p, -1, +1)$ after appending the three framed unknots to the chain. We obtain $A_{m+1} = -p\,A_{m} + A_{m-1}$, $A_{m+2} = -A_{m+1} + A_{m} = (p+1)\,A_{m} -A_{m-1}$ and $A_{m+3} = A_{m+2} + A_{m+1} = A_{m}$. So once again we start a cylic pattern $A_{m}, A_{m+1}, A_{m+2}, A_{m}$. The same holds for the $B_{m}$ and thus we obtain the following intersection number pattern for these framings: $-1, +1, -1$. If we change the orientation of $(A_m, B_m)$ we will obtain $+1, +1, +1$. $\Diamond$\\
\ \\
\noindent For example, $\frac{3}{2}$ $=(1,2) = [1,-2]$ forms the first element in a triad whose underlying curves are $[1,-2,-1] = \frac{2}{1}$ and $[1,-2,-1,-1] = 1$. The corresponding oriented filling curves are $(3,2)$, $(-2,-1)$ and $(-1,-1)$. Note that this argument gives an means for describing the cobordisms using relative Kirby calculus.

\section{Maps from triads}
\noindent Let $\LL$ be a link in a closed, oriented smooth three manifold $Y$. Suppose $\LL$ has $n$-components, ordered in some fashion, and suppose a triad is chosen in $\partial\big(Y - int\,N(L_{i})\big)$ for $ 1 \leq i \leq n$. Choose one curve from each triad, and label it $a$. This determines the labels on the other two curves. We will call $\{a,b,c\}^{n}$ the {\em code space}. Assign to each {\em code} $I \in  \{a, b, c\}^{n}$ a three manifold $Y(I)$ found by filling in the boundary of $Y - int\,N(\LL)$ using the framing determined by the code for each link component. \\
\ \\
\noindent An {\em immediate successor} to a code $I$ is a code $I'$ which agrees with $I$ in all but one spot and in that spot is one step larger in the ordering $a < b < c < a$. We can partially order $\{a, b, c\}^{n}$ using $a < b < c$ and extend by the lexicographic ordering on products. This also imposes an order on $\mathcal{C} = \{a,b\}^{n}$, an important subset of the codes. In fact, we can also order $\{c, a, b\} \times \{b, c, a\}$, or any other such product, by using the relevant three term inequalities from the cyclic ordering and starting from the first element of each factor. \\ 
\ \\
\noindent Let $\Gamma$ be a bouquet in $Y$ for the link underlying $\LL$. This is a one complex consisting of embedded, disjoint arcs connecting each component of $\LL$ to a specified basepoint in $Y$. We can always form a Heegaard diagram subordinate to $\Gamma$ where \\

\begin{defn}
Let $\Gamma$ be a bouquet for $\LL$. A pointed Heegaard diagram $\Siab$ is subordinate to $\Gamma$ if
\ben
\item $\Siab$  with $\be_{1}, \ldots, \be_{n}$  deleted describes $Y - \mathrm{int}N(\Gamma)$
\item After surgering $\be_{n+1}, \ldots, \be_{g}$, $\be_{i}$ lies in $F_{i}$, a punctured torus in $\partial N(\Gamma)$ surrounding $L_{i}$ (for each $i= 1, \ldots, n)$
\item $\be_{i}$ is a meridian for $L_{i}$
\een
\end{defn} 
\ \\
\noindent We delete $\be_{1}, \ldots, \be_{n}$ and replace them so that $\be_{i}$ represents the framing curve associated with $F_{a}^{i}$, the $a$-framing, and sitting in the punctured torus neighborhood in the previous definition. We let $\ga_{i}$ be a small Hamiltonian isotope of $\be_{i}$ for $i > n$ and a curve representing the framing $F_{b}^{i}$ for $i \leq n$, and finally we let $\delta_{i}$ be a curve representing a small Hamiltonian isotopes of $\be_{i}$ for $i > n$ and the framing $F_{c}^{i}$ for $i \leq n$. We will assume that all attaching circles representing different framings from a triad are chosen to intersect in one point {\em geometrically}. Finally, we assume that $F_{c}^{i}$ can be deformed to $F_{a}^{i} \cup F_{b}^{i}$ inside our Heegaard surface $\Si$. To see that this always possible, start with a Morse function from $Y- N(\Gamma)$ to $(-\infty,1]$ with the boundary uniformly sent to $1$. We can then add the attaching circles for the handles as specified, using a small isotoped smoothing of $F_{a} \cup F_{b}$ to represent $F_{c}$. \\
\ \\
\noindent For each code $I \in \{a,b,c\}$ we can form a new set of attaching circles $\eta(I)$ where
$$
\eta_{i}(I) = \left\{ \begin{array}{cl} 
													\be_{i} & i > n \mathrm{\ or\ } I_{i} = a \\ 
													\ga_{i} & \mathrm{if\ } I_{i} = b\\ 
													\de_{i} & \mathrm{if\ } I_{i} = c
											\end{array} \right.
$$
Then $(\Si, \seta, \eta(I), w)$ will represent $Y(I)$. We need to be more specific when we consider these embedded in the same Heegaard diagram. First, we place a special point $w$ in the diagram. We do this so that we can join $w$ to a point on each of $\be_{1}, \ldots, \be_{n}$ by paths which do not intersect any $\al$ or $\be$ curve. This is always possible. Furthermore, when forming $\eta(I)$, we need copies of various curves. Each copy of a curve will, in fact, be a small Hamiltonian isotope of the curve, and should intersect each of the other copies in only two points. These also should not cross the paths chosen from $w$. Furthermore, these are chosen so that every diagram we consider is weakly admissible.

\begin{center}
\begin{figure}
\includegraphics[scale=0.7]{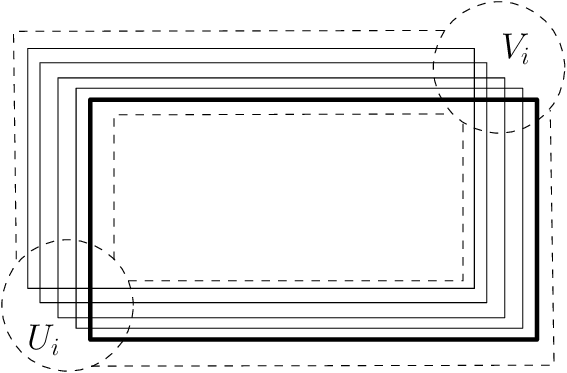}
\caption{A schematic representation of the annular neighborhood of one of the $\be_{i}$ curves (the thick curve) with various replicas of it, pushed off by small Hamiltonian isotopies. Notice that all the isotopes sit on one side or the other of the curve between $U_{i}$ and $V_{i}$, crossing only within those neighborhoods. The curve closest to $\be_{i}$, away from $U_{i}, V_{i}$, will always receive the lowest label, and the labels increase until we arrive at the curve farthest from $\be_{i}$.}\label{fig:specialpoints}
\end{figure}
\end{center}

\begin{center}
\begin{figure}
\includegraphics[scale=0.5]{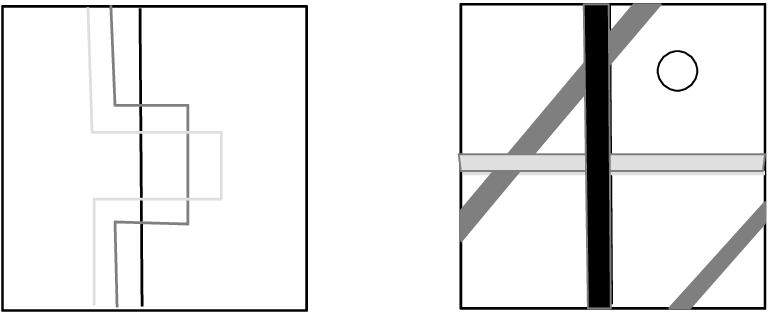}
\caption{We depict the conditions on the $\be$, $\ga$, and $\de$ curves, and their deformations. On the left, we show how we intend to isotope curves in a small annular region of each model curve. On the right, we depict neighborhoods of a triad of framings, thought of as $\infty, 0, -1$-framings. The torus is connected to the rest of the diagram inside the circle. There may be more than one connecting tube; however, the dark gray region can be deformed into the union of the other to regions. We place the isotopes of all the framing curves in the triad within the shaded annular neighborhoods.
Furthermore, the intersections of all the curves occur in a grid where two shaded neighborhoods overlap.}\label{fig:samplecurves}
\end{figure}
\end{center}

\noindent More precisely, we choose two points on each $\be_{i}$, and on each framing curve, and label them $1_{i}$ and $2_{i}$ (for the framing
curves these will also be labelled with $a,b$ or $c$). In a small annular region of the curve $\be_{i}$, small enough that all other curves cross
transversely, we choose two regions $U_{i}$ and $V_{i}$ such that $1_{i} \in U_{i}$ and $2_{i} \in V_{i}$. We then require that the copies of $\be_{i}$ are disjoint and to the right of $\be_{i}$ in the annular neighborhood, along one arc from $U_{i}$ to $V_{i}$ and disjoint and to the left between $V_{i}$ and $U_{i}$, i.e. along the other arc. Note that by isotoping $1_{i}$ to $2_{i}$ and $2_{i}$ to $1_{i}$ in the same direction, we can switch the roles of the two arcs, so it does not matter how we choose the arcs. Furthermore, we have asked that each pair of isotoped curves intersect in only two points, one necessarily in each of  $U_{i}$ and $V_{i}$, and we ask that intersections occur in order as depicted in the local diagram in Figure \ref{fig:specialpoints}. When we form the Hamiltonian isotopes of the framing curves, we require that all occur in neighborhoods of each curve that are sufficiently small that they reflect the same combinatorial structure on the surface as the original framing curves (i.e. there is a deformation retraction of the isotopes onto the union of the original framing curves). Since we have assumed that $\be_{i} \cap \ga_{i}$, $\ga_{i} \cap \de_{i}$ and $\de_{i} \cap \be_{i}$ consist of precisely one point, the isotopes of each type ($\be, \ga$ or $\de$) intersect any isotope of the other type for the same triad in precisely one point. All the isotopes of $\de$ will also be deformable into the annular neighborhood of $F_{a} \cup F_{b}$ for each link component. See Figure \ref{fig:samplecurves}.

\noindent Given a sequence of codes $\codeS{0}{k}$ we can construct a map
$$
D_{\codeS{0}{k}} : \hatc{Y(I^{0})} \lra \hatc{Y(I^{k})}
$$
This will only be a chain map when $k \leq 2$.\\
\ \\
\noindent Using the fact that $Y_{\eta(I^{j})\eta(I^{j+1})}$ is homeomorphic to a connect sum of $S^{1} \times S^{2}$'s, there is a generator of $HF^{\leq 0}$ with maximal absolute $\Q$-grading, which we will denote $\hth{j}{j+1}$. We then let
$$
 D_{\codeS{0}{k}}(\xi) = \widehat{f}_{\al, \eta(I^{0}), \cdots, \eta(I^{k})}(\xi \otimes \hth{0}{1} \otimes \cdots \otimes \hth{k-1}{k})
$$
\ \\
where the map on the right side is\\
$$
\widehat{f}_{\nu^{0}, \ldots, \nu^{m}}( {\bf x} \otimes \Theta_{0,1} \otimes \cdots \otimes \Theta_{k-1,k}) = 
\sum_{{\bf y} \in \T_{\nu^{0}} \cap \T_{\nu^{m}}} \sum_{\{\phi\,|\,\mu(\phi) = 0, n_{w}(\phi) = 0\}} \#\mathcal{M}(\phi){\bf y}
$$
Here the map is built out of counting pseudo-holomorphic $k+2$-gons in the homotopy class $\phi$ in $Sym^{g}(\Si)$ with boundary conditions specified by the intersection points above and the totally real tori $\T_{\al}, \T_{\eta(I^{0})}, \ldots, \T_{\eta(I^{k})}$ and Maslov index $0$. For our choice of curves $\Theta_{i,j}$ is always represented by a single generator in the {\em chain} complex, not just in homology. \\

\section{The $A_\infty$-relation in Heegaard-Floer homology}

\noindent The purpose of this section is to fix some terminology and motivate a result from 
the Floer theory of Lagrangian intersections. We do not present the details; for those, the
reader should consult \cite{Fooo} or \cite{Seid}.\\
\ \\
\noindent Let $\nu_{0}, \nu_{1}, \ldots, \nu_{k}$ be sets of attaching circles in $\Si$, i.e.
simple, disjoint curves which are linearly independent in $H_{1}(\Si)$. Let $\T_{i}$ be
the corresponding totally real torus in $Sym^{g}(\Si)$. We assume that the sets of attaching
circles have been chosen so that any sub-collection is weakly admissible in the sense that
any periodic domain which is the sum of doubly periodic domains has both positive and negative
multiplicities. Furthermore, we assume a generic choice of almost complex data: $J : \mathcal{M}(\triangle) \times \triangle \ra \mathcal{U}$.
Here $\triangle$ stands for any polygon, $\mathcal{M}(\triangle)$ is the moduli space of conformal structures on the polygon, moduli
conformal reparametrization, and $\mathcal{U}$ is a special class of almost complex structures defined in \cite{Hom3}, for a given
conformal structure $\mathfrak{j}$ on $\Si$.\\
\ \\
\noindent We consider the $1$-dimensional moduli spaces representing a homotopy class in $\pi_{2}({\bf x}, {\bf u_{1}}, \ldots,$ $ {\bf u_{k-1}}, {\bf y})$ where $u_{i} \in \T_{i} \cap \T_{i+1}$ and ${\bf x} \in \T_{0} \cap \T_{1}$ and ${\bf y} \in \T_{0} \cap \T_{k}$. These moduli
spaces have compactifications, a point we take as given, and the boundaries of these compactifications include those which come from broken
polygons (there may be others, but these will all cancel in the end and will not be considered here; the only non-trivial contributions will
come from these broken polygons). To specify the terminology, a broken polygon will arise from a {\em division} along two marked
edged, $\nu_{i}$ and $\nu_{j}$ if the homotopy class it represents arises from a splicing of the form $\psi_{1} \ast \psi_{2}$ where
$\psi_{1} \in \pi_{2}(\nu_{0}, \ldots, \nu_{i}, \nu_{j}, \ldots, \nu_{k})$ and $\psi_{2} \in \pi_{2}(\nu_{i}, \ldots, \nu_{j})$. Pictorially,
if $\triangle$ is a $k+1$-gon, with $\nu_{0}$ on one edge oriented counterclockwise, and $\nu_1, \ldots, \nu_{k}$ labelling the other edges
clockwise from $\nu_{0}$, this corresponds to choosing a chord between the $\nu_{i}$ and $\nu_{j}$ edges and contracting it. 

\begin{center}
\begin{figure}
\includegraphics[scale=0.7]{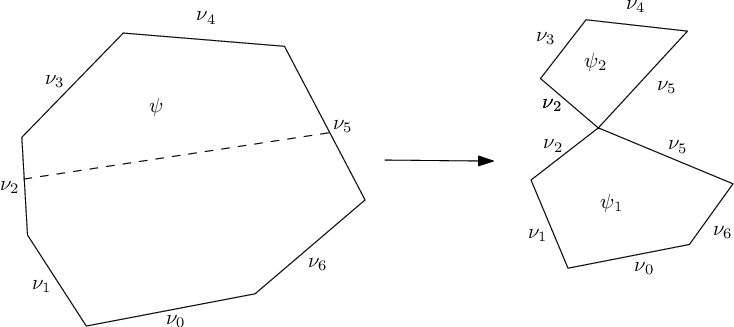}
\caption{Here $\psi = \psi_{1} \ast \psi_{2}$. The chord determines the edge labels for the homotopy classes which splice to $\psi$.}\label{fig:division}
\end{figure}
\end{center}

\noindent All the broken polygons in the compactification of a $1$-dimensional space have boundary contributions from a single division in 
the polygon, since the resulting parts will each have dimension $0$. Any further divisions results in pieces with negative formal dimension, and
thus do not have holomorphic representatives. Thus, $\#\mathcal{M}(\psi) = \#\mathcal{M}(\psi_{1}) \times \#\mathcal{M}(\psi_{2})$, where we are
counting components (over $\F_{2}$!). We note that divisions between consecutive edges of the polygon give rise to digons which contribute to
the differential of the Heegaard-Floer homology of the pair. We will assume, as will be the case in this paper, that the ${\bf u}_{i}= \hth{i}{i+1}$ in their respective Heegaard-Floer homologies, $HF^{\leq 0}(\T_{i}, \T_{i+1})$. We fix ${\bf x}$, and allow ${\bf y}$ to vary over all the possible intersection points for its pair. If we consider the $1$-dimensional moduli spaces for all homotopy classes which arise under these assumptions, and form the chain maps as in the preceding section, we arrive at the $A_{\infty}$ relation:
$$
\sum_{0 \leq i \leq j \leq k} \widehat{f}_{\eta(I^{i}), \ldots, \eta(I^{j})} \circ \widehat{f}_{\eta(I^{0}), \ldots, \eta(I^{i}), \eta(I^{j}), \ldots, \eta(I^{k})} = 0
$$

\section{Defining the chain complex}

\noindent We can use the maps $D_{\codeS{0}{1}}$ and the $A_{\infty}$-relation to define several chain complexes. \\

\begin{defn}
A {\em complete} subset of the code space, $\mathcal{S}$, is any subset such that for $I,J \in \mathcal{S}$, the following is always true $$\big(I < K < J \big) \Rightarrow \big(K \in \mathcal{S}\big)$$   
\end{defn}
\ \\
\noindent We note that both $\{a,b,c\}^{n}$ and $\mathcal{C}$ are complete subsets of the code space. \\ 
\ \\
\noindent For any complete subset of the code space, $\mathcal{S}$, let $X(\mathcal{S}) = \bigoplus_{I \in \mathcal{S}} \hatc{Y(I)}$. We will equip 
$X(\mathcal{S})$ with a map $D_{S}$, which will be our differential, by defining for $\xi \in \hatc{Y(I)}, I \in\mathcal{S}$
$$
D_{\mathcal{S}}(\xi) = \sum_{J \in \mathcal{S}}\  \sum_{I = \codeS{0}{k} = J} D_{\codeS{0}{k}}(\xi)
$$
where the sequence $\codeS{0}{k}$ consists of an increasing sequence of immediate successors. The completeness property
is necessary to ensure that such sequences can only occur using codes in $\mathcal{S}$. \\

\begin{thm}
$D^{2} \equiv 0$ on $X(\mathcal{S})$ for any complete subset, $\mathcal{S}$, of $\{a, b, c\}^{n}$.
\end{thm}
\ \\
\noindent {\bf Proof:} Let $\xi \in \hatc{Y(I^{0})}$. We wish to compute $D^{2}(\xi)$. We can simplify matters by instead calculating the coefficient of $\nu$ in $\hatc{Y(I^{k})}$ for $I^{k} \geq I^{0}$. This occurs in the composition, $D^{2}$, only for sequences of immediate succesors with
$I^{0} < I^{1} < \cdots < J = I^{j} < \cdots < I^{k}$ where the first application of $D$ provides a sum over sequences of the form $I^{0} < \cdots < J$ and the second application provides a sum for $J < \cdots < I^{k}$. Since $I^{0}$ and $I^{k}$ are both in $\mathcal{S}$, all such sequences consist of elements of $\mathcal{S}$ and thus appear in the sums defining $D$. The coefficients that appear in $D^{2}$ are sums of products $\#\mathcal{M}(\phi_{1})\#\mathcal{M}(\phi_{2})$ where the moduli spaces are for pseudo-holomorphic polygons, the first of which originates at $\xi$, goes through the canonical generators of the $Y_{\eta(I^{j}I^{j+1})}$, and terminates at $\xi'$ for some $\xi' \in \hatc{Y(J)}$, while the second would then start at $\xi'$ and terminate at $\nu$. From the gluing theory in \cite{Seid} these are the boundaries of one dimensional moduli spaces for $\pi_{2}(\xi, \hth{0}{1}, \cdots, \hth{k-1}{k}, \nu)$. Thus, these occur as the compositions in the $A_{\infty}$ relation. For a given immediate successor sequence $\codeS{0}{k}$, the $A_{\infty}$ relation provides:
$$
\sum_{-1 \leq i \leq j \leq k} \widehat{f}_{\eta(I^{i}), \ldots, \eta(I^{j})} \circ \widehat{f}_{\eta(I^{-1}), \ldots, \eta(I^{i-1}), \eta(I^{i}), \eta(I^{j}), \ldots, \eta(I^{k})} = 0
$$
where $\eta(I^{-1}) = \alpha$. The terms which occur in $D^{2}$ are precisely those where $i = -1$. To account for the other terms, observe that when $i \neq -1$ there can be many paths of immediate successors $\codeS{i}{j}$. In $D^{2}$ we would sum over all these
different paths. It will suffice to show that, when $i \neq -1$, the total is zero in this sum. The sum of all the $i=-1$ terms provide all the terms in $D^{2}$, so the identity above will show that $D^{2} \equiv 0$. We postpone the remainder of the proof to accumulate some useful background results.\\

\subsection{Background}

\begin{lemma}\label{lem:workhorse}
Let $\Siabg$ have $\Si_{\be\ga} \cong \#^{l} S^{1} \times S^{2}$. Pick classes $[x] \in \widehat{HF}(\Si_{\al\be})$ and $[y] \in \widehat{HF}(\Si_{\al\ga})$ representing torsion $Spin^{c}$ structures on the boundary and which are $\Q$-grading homogenous. 
Suppose that $\langle \widehat{F}_{\Sp{u}}([x]), [y] \rangle \neq 0$. Then for any generator, ${\bf x}$ of $\widehat{CF}(\Si_{\al\be})$ which occurs with non-zero coefficient in a linear combination representing $[x]$, and for any generator ${\bf y}$ with a similar property for $[y]$,
we have that $\mu(\psi) = 0$ for any homotopy class in $\pi_{2}({\bf x}, \Theta, {\bf y})$ representing $\Sp{u}$ and having $n_{w}(\psi) = 0$.
\end{lemma}
\ \\
\noindent {\bf Proof:} According to the assumption about the cobordism map, there is at least one pair of intersection points ${\bf x}'$, and ${\bf y}'$ representing their respective $Spin^{c}$ structures, and one homotopy class $\psi$ joining them, for which the conclusion holds. Any other pair of intersection points can be joined to these by paths in the boundary. However, since each of these generators descends to a homology class with homogenous $\Q$-grading, each of the generatorsmust also have that $\Q$-grading. Thus they can be joined by a path with both $n_{w} = 0$ and $\mu = 0$. As a result, there is at least one homotopyclass of triangles joining ${\bf x}$ and ${\bf y}$, representing $\Sp{u}$ and having $n_{w} = 0$ and $\mu = 0$. All other homotopy classes of triangles with this property differ by doubly periodic domains in the boundary. However, since the $Spin^{c}$ structures on the boundary are all torsion, the Maslov indices of these doubly periodic domains are all $0$. Thus the conclusion holds for all the relevant homotopy classes of triangles. $\Diamond$. \\
\ \\
\noindent We can apply this lemma to the following types of triangles:\\
\ben
\item Where all three boundaries are either $S^{3}$'s or connected sums of $S^{1} \times S^{2}$'s and the cobordism represented is $\big(\#^{l} S^{1} \times S^{2}\big) \times I$. For example, these triples occur in the proof of the invariance of Heegaard-Floer homology under handleslides of the attaching curves.
\item Where all three boundaries are as in the previous listing, but the cobordism represents surgery on a homologically non-trivial, but primitive, 
knot in $\Si_{\al\be}$.
\item Where all three boundaries are the same as before, and the cobordism represents $-1$ surgery on an topologically trivial unknot inside $\#^{l} S^{1} \times S^{2}$, i.e. a blow up of the trivial cobordism. 
\een
\ \\
\noindent That the previous lemma applies to these are standard results in Heegaard-Floer homology. We make a few comments, however. For the first two listings, there is only one $Spin^{c}$ structure joining the non-trivial homology generators on the boundary. This structure is the torsion one, since all of $H_{2}$ arises in the boundary. That the map on homology is non-trivial in the first case comes from invariance preserving $\Q$-grading. In the second, it comes from Prop. 9.3 of \cite{AbsG}:
\begin{lemma}\cite{AbsG}
Let $Y$ be a closed oriented three-manifold and $K \subset Y$ be a framed knot. Assume that the cobordism $W_{K}$, found by adding a two handle to $Y\times I$ along $K$, has $b_{2}^{-}(W_{K}) = b_{2}^{+}(W_{K})$. Let $\Sp{s}$ be a $Spin^{c}$ structure on $W_{K}$ whose restriction, $\Sp{t}$, to $Y$ and its restriction, $\Sp{k}$, to $Y_{K}$ are both torsion. When $K$ represents a non-torsion class in $H_{1}(Y;\Z)$, then if $HF^{\infty}(Y,\Sp{t})$ is standard, the induced map
$$
F^{\infty}_{W_{K},\Sp{s}}: HF^{\infty}(Y,\Sp{t}) \lra HF^{\infty}(Y_{K},\Sp{k})
$$
vanishes on the kernel of the action of $[K] \in H_{1}(Y;\Z)$ and induces an isomorphism
$$
HF^{\infty}(Y,\Sp{t})/\mathrm{Ker}[K] \cong HF^{\infty}(Y_{K}, \Sp{k})
$$
If $K$ represents a torsion class in $H_{1}(Y;\Z)$ and $HF^{\infty}(Y_{K},\Sp{k})$ is standard, then the map
$F^{\infty}_{W_{K},\Sp{s}}$ induces an isomorphism 
$$
HF^{\infty}(Y,\Sp{t}) \cong \mathrm{Ker}\big( [L] : HF^{\infty}(Y_{K}, \Sp{k}) \ra HF^{\infty}(Y_{K}, \Sp{k})\big)
$$
where $[L] \in H_{1}(Y_{K},\Z)$ is the core of the glued-in solid torus.
\end{lemma}
\noindent To complete the argument, note that the grading change formula forces the same conclusion to hold on $\widehat{HF}$ in the cases we consider.\\
\ \\
For the last, there is a homology class in $H_{2}$ for the cobordism not arising from the boundary. However, the map is fully understood through the blow-up formula, \cite{Smoo}, which in turn rests on the following lemma:
\begin{lemma} \cite{Smoo}
Let $(E, \al_{0}, \be_{0}, \ga_{0})$ be a genus $1$ Heegaard triple representing the cobordism between $S^{3}$'s found by $-1$ surgery on an unknot in $S^{3} \times I$. Thus, $[\ga_{0}] = [\be_{0}] - [\al_{0}]$. Assume that the attaching circles are in general position, pairwise intersecting once geometrically. For each non-negative integer $k$ there are two homotopy classes of triangles $\psi^{\pm}_{k}$ such that $\mathcal{M}(\psi_{k}^{\pm})$ consists of a single point, $n_{w}(\psi_{k}^{\pm}) = \frac{k(k+1)}{2}$ and $\langle c_{1}(\Sps{s}{w}(\psi_{k}^{\pm}), [\mathcal{T}] \rangle = \pm (2 k + 1)$.
\end{lemma}
\noindent Namely, the blow-up map for each $Spin^c$ structure on $Y\times I \# \mir{{\C}P}^{2}$ is multiplication by $U^{i(i+1)/2}$ where $i$ is determined by the pairing of the first Chern class of the $Spin^{c}$ structure with the $-1$-sphere: $\langle c_{1}(\Sp{u}), S \rangle - 1 = 2i$. The map is trivial for the $\hat{\ }$-theory except when $i = 0, -1$. In each of these cases, lemma \ref{lem:workhorse} applies. More importantly, only these cases can arise in $\widehat{F}$. \\
\ \\
\noindent Finally, we include a gluing result which is necessary for decomposing our diagrams into diagrams that are easier to analyze:
\begin{lemma}
Fix a pair of Heegaard diagrams $(\Si, \eta_{0}, \ldots, \eta_{k})$ and $(E, \nu_{0}, \ldots, \nu_{k})$ where the latter is genus $1$. Form the Heegaard diagram $\Si \# E$ with the attaching circles as above. Let $\sigma \in \Si \# E$ be a point in the connect sum region. Let $\psi$ be a homotopy class of $k+1$-gons on $(\Si, \eta_{0}, \ldots, \eta_{k})$ and $\psi_{0}$ be one on $(E, \nu_{0}, \ldots, \nu_{k})$. Assume that 
$\mu(\psi) = 0$, $\mu(\psi_{0}) = k - 2$, and $n_{\sigma}(\psi_{0}) = 0$. Furthermore, assume that $\mathcal{M}(\psi_{0})$ is isomorphic to copies of
the moduli space of conformal structures on a $k+1$-gon in the complex plane. Let $\psi'$ be the homotopy class of $k+1$-gons for the diagram $\Si \# E$ with the multiplicities of $\psi$ on $\Si$ and $\psi_{0} + n_{\sigma}(\psi)[E]$ on $E$. Then for a suitable choice of families of complex structures and perturbations there is an isomorphism
$$
\mathcal{M}(\psi') \cong \mathcal{M}(\psi) \times \#\big(\mathcal{M}(\psi_{0})\big)
$$
where $\#\mathcal{M}(\psi_{0})$ is the set of components to the moduli space of $k+1$-gons for $E$.
\end{lemma}

\noindent The point is that we can match any conformal structure of $k+1$-gon in the moduli space of $\psi_{0}$. For each copy of this moduli space in that of $\psi_{0}$ we can obtain a match with a choice of representative in $\mathcal{M}(\psi)$.   

\subsection{Finishing $D^{2} \equiv 0$}\ \\
\ \\
\noindent To complete the proof of $D^{2} \equiv 0$ we only need the following lemma:

\begin{lemma} Fix $I, J \in \{a,b,c\}^{n}$. Then 
$$
\sum_{I = \codeS{0}{k}=J} \widehat{f}_{\eta(I^{0})\cdots\eta(I^{k})}(\hth{0}{1} \otimes \cdots \otimes \hth{k-1}{k}) = 0
$$
where the sum is over all immediate successor sequences between $I$ and $J$.
\end{lemma}

\noindent {\bf Proof:} First suppose $k > 2$. We then consider the triples $(\eta(I^{0}), \eta(I^{j}), \eta(I^{j+1}))$. These represent cobordisms that are either surgery on a non-trivial primitive homology class in $\#^{l} S^{1} \times S^{2}$ or blow-ups of a trivial cobordism. The former occur when the alteration from $I^{j}$ to $I^{j+1}$ occurs on a curve which in $I^{j}$ is a small Hamiltonian isotope of one in $I^{0}$, i.e. an $a \ra b$ change in code from from $I^{j}$ to $I^{j+1}$ or a change  $b \ra c$ from $I^{0}$ and $I^{j}$ to $I^{j+1}$. The blow-up occurs when there is a change from $I^{0}$ to $I^{j}$ on the same indexed curve which changes in going to $I^{j+1}$. Furthermore, each of the boundaries of these triples consist of connect sums of some number of copies of $S^{1} \times S^{2}$. The lemma in the previous section implies that any homotopy class relevant to the calculation of the map $\langle \widehat{F}(\hth{0}{j} \otimes \hth{j}{j+1}) , \hth{0}{j+1} \rangle$ must have Maslov index zero. Furthermore, as above, there are relevant homotopy classes with non-trivial moduli spaces, and $\widehat{F}(\hth{0}{j} \otimes \hth{j}{j+1}) = \pm \hth{0}{j+1}$ for 
the surgery on the non-trivial homology class, while $\widehat{F}$ equals $0$ (see below) for the $-1$-sphere. We can construct a homotopy class for $\pi_{2}(\hth{0}{1}, \cdots, \hth{k-1}{k}, \hth{0}{k})$ by splicing together homotopy classes of triangles from these triples. This homotopy class would necessarily have $\mu = k - 2 > 0$ from the increased dimension from splicing. We can alter this homotopy class while preserving $n_{w} = 0$ only by adding doubly periodic domains from $\Si_{\eta(I^{j})\eta(I^{m})}$, which do not change the Maslov index since the $Spin^{c}$ structure will be torsion on these, or by adding triply periodic domains from the $-1$-surgery spheres. The latter change the first Chern class pairing by multiples of $2$, a difficulty we address forthwith. Finally, since the other generators of $\widehat{HF}(\#^{l} S^{1} \times S^{2})$ have lower grading, 
any homotopy class in $\pi_{2}(\hth{0}{1}, \cdots, \hth{k}{k+1} , \Theta^{-}_{0,j+1})$ will have {\em larger} Maslov index. For examples of spliced
homotopy classes see Figure \ref{fig:splice}.

\begin{center}
\begin{figure}
\includegraphics[scale=1.0]{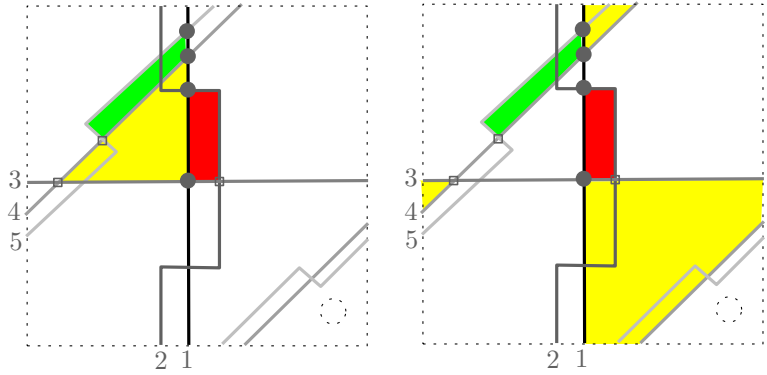}
\caption{We depict two homotopy classes found from splicing in a genus $1$ component including all three framing curves. The dotted circle delineates the connect sum region for the component. The three homotopy triangles have domains indicated in the yellow, green, and red. The dots are the intersection points along the $\eta(I^{0})$ curves.}\label{fig:splice}
\end{figure}
\begin{figure}
\includegraphics[scale=1.0]{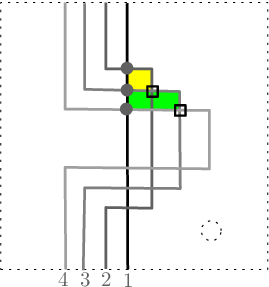}
\caption{Another example of a splicing of triangles. This type of configuration occurs in all the components coming from non-framing $\be$'s}\label{fig:splice2}
\end{figure}
\end{center}
\ \\
\noindent However, there is a class with $\mu = 0$ representing any $Spin^{c}$ structure on a triple where the variation in the attaching curves occurs all in one triad, and thus includes all three triad curves. This is the situation for the cobordism map for $-1$ surgery on the unknot. These have $n_{w} = \frac{i(i+1)}{2}$ for some $i$. As we have seen, only $i = 0, -1$ contribute to the cobordism map for the $\hat{\ }$-theory. However, by subtracting copies of the Heegaard surface, we can get classes with $n_{w} = 0$ and $\mu = -i(i+1)$. The argument above might fail since it depends on the Maslov index increasing under splicing, whereas here the formal dimension of the moduli space could decrease after splicing. However, for the choice of framing curves above we need not worry about this in the $\hat{\ }$-theory. Since $n_{w}(\phi) = 0$ in the whole $k+1$-gon, the non-zero multiplicities must occur in genus $1$ summands of the Heegaard multituple diagram. Lifting to the universal cover, and recalling that we can deform to the original framing curves, shows that the only possibilities have multiplicities $-i$ and $1+i$ in the two large triangles formed between the annular neighbrhoods of the framing curves, see Figure \ref{fig:samplecurves}. Of course, near the framing curves the Hamiltonian isotopes will change the multiplicities, but to avoid $w$, these are the only multiplicities for the triangles that are possible. Thus, $i =0, -1$ since we have non-negative multiplicity. We already know how to construct homotopy classes with this property by splicing the canonical triangles. Any other homotopy class through the same intersection points differs either in $i$, or by copies of the Heegaard surface (changing $n_{w}$), or by doubly periodic domains. Since we started with the canonical small triangles, adding doubly periodic domains will introduce negative multiplicities. Thus there are only two possible homotopy classes with $n_{w} = 0$ and $\mathcal{D} \geq 0$ in this genus $1$ summand, and both have Maslov index $k - 2$. Since this is true in all of the genus $1$ summmands, and $k-2$ is the dimension of the moduli space of conformal structures on a $k+1$-gon, we have that the overall Maslov index is also $k-2$. \\
\ \\
\noindent Hence we need only consider those classes with $k \leq 2$. For $k=1$, the proposition follows from $\hth{j}{j+1}$ being a cycle. While for $k=2$, there are two cases to consider. In the first, the terms in $I = I^{0} < I^{1} < I^{2} = J$ differ in only one framing curve. This is the $-1$-sphere case, and it is known that the triangles for $i = 0$ and $i = -1$ both have a single element in their moduli spaces, and these elements
cancel. In the other case $I$ and $J$ differ in two places. So there are exactly two codes $K$ and $L$ which lie between $I$ and $J$ depending upon which place we change first. We label these so that $K < L$ then, over $\F_{2}$, 
$$
\begin{array}{c}
\widehat{f}_{\eta(I)\eta(K)\eta(J)}(\hth{I}{K} \otimes \hth{K}{J}) = \hth{I}{J}\\
\ \\
\widehat{f}_{\eta(I)\eta(L)\eta(J)}(\hth{I}{K} \otimes \hth{K}{J}) = \hth{I}{J}\\
\end{array}
$$
The maps are as described since the diagrams for these triples reduce to genus 1 summands, all but two of which are the identity. In the other two
we will have a small Hamiltonian isotope followed by a change in framing, or vice-versa. In either case there is a unique triangle class represented
by a single element moduli space joining the canonical generators. $\Diamond$ 
\begin{center}
\begin{figure}
\includegraphics[scale=1.0]{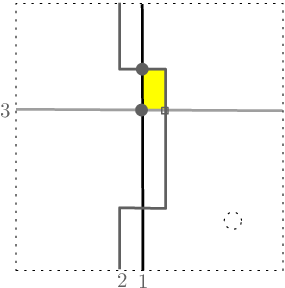} \hspace{.5in} 
\includegraphics[scale=1.0]{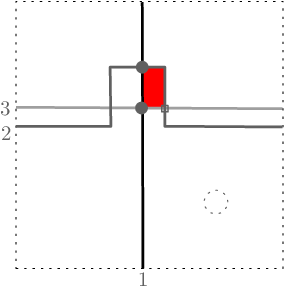}
\caption{The local pictures of a genus $1$ component for a framing which changes in $I < J < K$, and changes occur in two different triads. The left side corresponds to a code change of $00 \ra 01 \ra 11$ where the component depicted corresponds to the first label. The second diagram corresponds to $00 \ra 10 \ra 11$. The component corresponding to the second entry has the diagrams reversed. Notice that the holomorphic triangles have $1-3$ vertex at the same intersection point.}\label{fig:TypeIII}
\end{figure}
\end{center}
\noindent {\bf Comment:} We note that the proof for $X(\mathcal{C})$ is substantially similar, since we only allow two of the three framings in the triad. Thus, there are no triply periodic domains in any of the triangles, and there is no need to consider the homotopy classes of triangles
representing different $spin^{c}$-structures. 
\begin{center}
\begin{figure}
\includegraphics[scale=1.0]{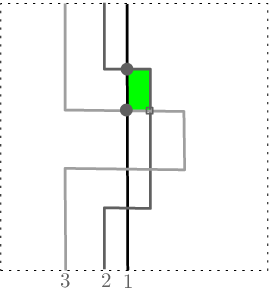}
\caption{For all the non-framing $\be$-curves, or any $\be$-curves whose code does not change, the genus $1$ diagram is as depicted, with a single unique holomorphic triangle. These components can be removed by the gluing lemma.}\label{fig:TypeII}
\end{figure}
\end{center}

\section{The spectral sequence}

\noindent Let $\mathbb{L}$ have an $n$-components each equipped with a triad, and let $\mathcal{S} = \{a,b\}^{n-1} \times \{a,b,c\}$ and let $\mathcal{C} = \{a,b\}^{n}$. We will show that

\begin{prop}
$$ H_{\ast}(X(\mathcal{S})) \cong 0 $$
\end{prop}

\noindent {\bf Proof:} There is  a short exact sequence of complexes
$$
0 \lra X(\{a,b\}^{n-1} \times \{c\}) \stackrel{i}{\lra} X(\mathcal{S}) \stackrel{\pi}{\lra} X(\mathcal{C}) \lra 0
$$
We will show that $H_{\ast}(X(\mathcal{S})) \equiv 0$ by proving that the connecting homomorphism is an isomorphism. The proof
relies upon the following homological lemma for coefficients in $\F_{2}$, whose proof is in the appendix:\\

\begin{lemma}\cite{Doub}
Let $\{A_{i}\}_{i=0}^{\infty}$ be a set of chain complexes and let $\{f_{i} : A_{i} \ra A_{i+1}\}_{i=0}^{\infty}$ be a set of chain maps
satisfying the properties:
\ben
\item For each $i$, $f_{i+1}\circ f_{i}$ is chain homotopically trivial, with chain homotopy $H_{i} : A_{i} \ra A_{i+2}$.
\item For each $i$, the map $f_{i+2}\circ H_{i} + H_{i+1} \circ f_{i} : A_{i} \ra A_{i+3}$ is a quasi-isomorphism. 
\een
then $M(f_{1})$ is quasi-isomorphic to $A_{3}$. 
\end{lemma}
\ \\
\noindent To prove $H_{\ast}(X(\mathcal{S})) \equiv 0$ we will use this lemma as follows. Let $\mathcal{G}$ be a complex, with chain group $A_{1} \oplus A_{2} \oplus A_{3}$ and with differential
$$
\partial = \left( \begin{array}{ccc} 
			\partial_{1} & 0 & 0 \\
			f_{1} & \partial_{2} & 0 \\
			H_{1} & f_{2} & \partial_3 \\
\end{array} \right)
$$
The first assumption about $\{A_{i}\}$ ensures that this is a differential.  This complex has homology $H_{\ast}(\mathcal{G}) \equiv 0$ since there is a short exact sequence $$ 0 \lra A_{3} \stackrel{i}{\lra} \mathcal{G} \lra M(f_{1}) \lra 0$$ yielding a long exact sequence with connecting homomorphism $H_{\ast}(M(f_{1})) \lra H_{\ast}(A_{3})$. This connecting homomorphism is precisely the quasi-isomorphism in the lemma. Thus, $H_{\ast}(\mathcal{G}) \equiv 0$. \\ 
\ \\
\noindent For the complexes arising from surgeries on $\mathbb{L}$, we will use $A_{1} \cong X(\{a,b\}^{n-1} \times \{a\})$, $A_{2} \cong X(\{a,b\}^{n-1} \times \{b\})$, and $A_{3} \cong X(\{a,b\}^{n-1} \times \{c\})$. For $i \geq 3$ we then repeat isomorphic copies of these cyclically. Let $\xi \in Y(I)$ for some $I$ with $a$ as the last entry. We define the chain map $F_{a}: A_{1} \ra A_{2}$ by
$$
F_{a}(\xi) = \sum_{J \in \{a,b\}^{n-1}\times \{b\}}\  \sum_{I = \codeS{0}{k} = J} D_{\codeS{0}{k}}(\xi)
$$
That this is a chain map follows from $D^{2} \equiv 0$ on $X(\mathcal{C})$. Namely if we let $\widetilde{D}_{i}$ be the differential on
$X(\{a,b\}^{n-1}\times \{i\})$, then the differential on $X(\mathcal{C}) \cong A_{1} \oplus A_{2}$ is 
$$
\left(\begin{array}{cc}
\widetilde{D}_{a} & 0 \\
F_{a} & \widetilde{D}_{b} \\
\end{array}\right)
$$
This realizes $X(\mathcal{C})$ as $MC(F_{a})$. We can likewise define $F_{b}$ on those groups with codes ending in $b$. Finally, we can describe $F_{c} : A_{3} \ra A_{4}$. To do this, we note that $A_{4} \cong A_{1}$ and so we can find a representative Heegaard multituple, which we choose to consist of small Hamiltonian isotopes of the attaching circles for $A_{1}$. Since each is really a direct sum of groups for different diagrams, the isotopies occur in these different diagrams. We do this for each $i \geq 4$ as well.  We can then define $F_{i}$ for $i \geq 3$ since the map has the same definition as $F_{a}$; we simply relabel which framing we choose to be $a$. In other words, the development until now has been symmetric in the three framings. We broke the symmetry by calling one of them $a$, but the theory works no matter which we call $a$. \\
\ \\
\noindent To complete the proof we must now verify that the conditions in the lemma hold. \noindent We can now verify that the compositions $F_{i+1} \circ F_{i}$ are all chain homotopic to $0$ and identify the homotopies. We once again consider compactness properties for $\mu =1$ spaces of pseudo-holomorphic $l$-gons, requiring that the last entry of the code for the target be two spots different from the source entry in the cyclic order $a < b < c < a$.  The compactness result above will apply to individual sets of generators, and using the same lemma we can cancel many terms -- those which do not correspond to $i = -1$; i.e. splicing homotopy classes both involving the $\alpha$'s -- by summing over the immediate successor sequences. This works because an alteration of the code by two spots allows us to still consider
$\{a, b, c\}^{n}$ as an $n$-dimensional cube (possibly with relabelling); we do not need to use the full cyclic ordering, yet. The terms with $\alpha$'s in each map in the composition fall into two types: 1) those where the last term in the code changes in each map, and 2) those where the last
element in the code alters by two spots in only one of the factor maps in the composition. The former correspond to the moduli spaces coming from 
$F_{i+1} \circ F_{i}$. The latter from compostions of the form $H_{a} \circ \widetilde{D}_{a}$ or $\widetilde{D}_{c} \circ H_{a}$ where $H_{a} : A_{1} \ra A_{3}$ can be written
$$
H_{a}(\xi) = \sum_{J \in \{a,b\}^{n-1}\times \{c\}}\  \sum_{I = \codeS{0}{k} = J} D_{\codeS{0}{k}}(\xi)
$$
By the standard argument $F_{b} \circ F_{a} = H_{a} \circ \widetilde{D}_{a} + \widetilde{D}_{c} \circ H_{a}$. By the cyclic symmetry in the setup,
this argument immediately implies that $F_{i+1} \circ F_{i}$ is chain homotopic to $0$ through a chain homotopy $H_{i}$ defined by looking at a two step change in the code, with maps defined relative to the Heegaard tuples described above.\\
\ \\
\noindent We need now to verify hypothesis (2): $F_{c} \circ H_{a} + H_{b} \circ F_{a}$ is a quasi-isomorphism. Again the argument we will use can be extended cyclically to establish the same conclusion for $F_{3} \circ H_{1} + H_{2} \circ F_{1}$. We note that the maps $H$ is found by counting pseudo-holomorphic $n$-gons which as we go around the boundary alter some framing curve by two steps in  $a \ra b \ra c \ra a$ depending upon the subscript for the map. To obtain the result we consider moduli spaces of pseudo-holomorphic $n$-gons which, as we go from $\al \eta(I^{0})$ to $\al\eta(I^{k})$ have, for some framing curve, the entire cycle $a \ra b \ra c \ra a$ in their boundary, again starting at a code determined by the map's subscript. Now we cannot mindlessly use the cancellation lemma above, since we cannot restrict attention to a lexicographically ordered subset of the code space: the cycle prevents that. Instead, we consider what maps arise in that $A_\infty$-identity. First, if we divide along the $\al$-edge then we have two possibilities: 1) all the alterations on that framing curve occur on one side of the division: the other side therefore comes from $\widetilde{D}_{a}$ or $\widetilde{D}_{4} \cong \widetilde{D}_{a}$, while the first side contributes to a map $G$ found by studying pseudo-holomorphic $n$-gons with $\mu =0$ and built using the entire cycle of framings (i.e. a three step change in a single code entry); 2) the division divides the steps in a code entry into $1$ and $2$ or $2$ and $1$, these give $F_{c} \circ H_{a} + H_{b} \circ F_{a}$. If the division occurs along two $\eta$-edges ($i \neq -1$ in the $A_{\infty}$-relation) there are also two possibilities: 1) the all $\eta$ homotopy class of polygons includes at most two of the three alterations, or 2) all three steps occur in the all $\eta$ homotopy class. In the former case the cancellation lemma still applies since we can still consider the immediate successor sequences in a lexicographically ordered subset of the code space. When we add over immediate successors all the terms cancel, implying that this situation contributes nothing. It is to the second case that we now must turn. 
\begin{center}
\begin{figure}
\includegraphics[scale=1.0]{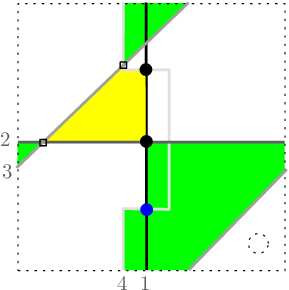} \hspace{0.5in} \includegraphics[scale=1.0]{TypeI}
\caption{We depict the genus $1$ diagram for the $a \ra b \ra c \ra a$ cycle in framing curves. The $a$ framing is the darkest, follwed by $b$, $c$, and $a'$ as the lightest. The dark circles are the generators for $ab$ and $a'b$. The squares are the canonical generators for $bc$ and $ca'$. Note that for the $bca'$ triple, the obvious holomorphic triangle has vertex at $\Theta^{-}_{ca'}$, the blue dot, which gives rise to a $\mu=1$ holomorphic quadrilateral, in green. However, there is a unique quadrilateral, in yellow, in the correct homotopy class, which allows us to complete the argument for $k=3$. Every other annular region consists of small Hamiltonian isotopes. The shaded quadrilateral found from splicing has $\mu = 1$ and is holomorphic. Moreover, by cutting along the obvious curves we can obtain any conformal structure in the one dimensional moduli space of quadrilaterals in $\C$. Using the gluing theorem above, we obtain a unique quadrilateral map.}\label{fig:hatquad}
\end{figure}
\end{center}
\noindent We can try repeating the argument for the cancellation lemma; however, one of the Heegaard triangles -- the one for the $c \ra a$ alteration -- does not have a pseudo-holomorphic triangle to use. Instead, Figure \ref{fig:hatquad} (see also \cite{Doub}) indicates that the appropriate triangle goes to $\Theta^{-}$. The triangle to $\Theta^{+}$, therefore, has $\mu = -1$. Splicing the triangle for this step does not change the Maslov index of the homotopy polygon to which it is spliced. So for this case the Maslov index is $k - 3$ and we need to concern ourselves about the $k=3$ case as well as the $k \leq 2$ case. For $k = 3$ we have a homotopy quadrilateral where all the steps in the code entry appear in the four edges. By the local calculation in Figure \ref{fig:hatquad}, there is a unique pseudoholomorphic quadrilateral which yields the map in this case. This map has image $\Theta^{+}_{00'}$.  \\
\ \\
\noindent Thus the map on $E^{1}$ will be the identity on the graded homology, since this map will be found by the triangle map from $Y(I)$ to $Y(I)$ where each of the $\ga$'s is a small isotope of the $\be$'s. In particular, the map induced by the chain map will be a quasi-isomorphism since the upshot of the above argument is that $F_{c} \circ H_{a} + H_{b} \circ F_{a} = G \circ \widetilde{D}_{a} + \widetilde{D}_{4} \circ G + I + N$ where $N$ lowers the filtration index found by sending $a \ra 0$ and $b \ra 1$ and summing over the code. Since the top level of $I + N$ is a isomorphism of the filtered complexes, and the filtration is bounded below, this chain map will induce an isomorphism on the chain complexes. $\Diamond$\\
\ \\
\noindent As a consequence of the preceding calculation, we deduce:

\begin{thm}
$$H_{\ast}(X(\mathcal{C})) \cong \widehat{HF}\big(Y(\{cc\ldots c\})\big)$$
\end{thm}

\noindent {\bf Proof:} We proceed by induction on $n$. For our base case, when $n=0$ and all the components are filled with the framing determined by $c$, we obtain $\widehat{HF}\big(Y(\{c\ldots c\})\big)$, and the conclusion follows directly. We apply the result above to the link $\mathbb{L}'$ obtained by taking the first $k$ components of $\mathbb{L}$ in the three manifold obtained by $c$-surgery on the remaining $n-k$ link components. 
We apply the previous proposition, that $X(\mathcal{S})$ has trivial homology, where the codes now apply to the first $k$ link components. From this
we obtain that 
$$
H_{\ast}(X(\{a,b\}^{k}\times\{c\ldots c\})) \cong H_{\ast}(X(\{a,b\}^{k-1}\times\{c\ldots c\}))
$$
Proceeding with the induction proves the result. $\Diamond$\\
\ \\
\noindent The case for $n=1$ has indpendent interest:

\begin{thm} \cite{Doub}
Let $K$ be a knot in a three manifold $Y$ equipped with a labelled triad as above. Let $$ \widehat{f}: \hatc{Y_{a}} \ra \hatc{Y_{b}}$$ be the Heegaard triple chain map. Then $\hatc{Y_{c}}$ is quasi-isomorphic to the $M(\widehat{f})$. 
\end{thm}
\ \\
\noindent For $\mathbb{L} \subset Y$, we note that $X(\mathcal{C})$ is filtered by taking the number of $b$'s in the code
$I \in \{a,b\}^{n}$ defining $\widehat{HF}(Y(I))$. Using the standard Leray spectral sequence for this filtration, we obtain\\

\begin{thm}\cite{Doub}
There is a spectral sequence whose $E^{1}$-term is $\bigoplus_{I \in \mathcal{C}} \widehat{HF}(Y(I))$ with the properties that
\ben
\item The differentials respect the lexicographic ordering on $\mathcal{C}$. 
\item The differential $d_{1}$ on the $E^{1}$-page for the summand corresponding to $I$ is found by adding the cobordism maps $\widehat{F}_{I < I'}$ over the immediate successors, $I'$, to $I$.
\item The spectral sequence converges to $\widehat{HF}\big(Y(\{cc\ldots c\})\big)$ in finitely many steps. 
\een 
\end{thm}

\part*{II. Additional facts about the link surgery spectral sequence} \noindent In this part we expand the construction of the link surgery spectral sequence
to restore some of the $Spin^{c}$ information. Once we have explained how to adapt the proof to the setting of knot Floer homology and certain twisted coefficient systems, we begin to answer a question from the introduction to \cite{Doub}. Namely, we establish the invariance properties of the spectral sequence under changes of bouquet, and develop a theory of cobordism maps between link surgery spectral sequences. These sections form the
original part of this paper. \\

\section{Adjustments for knot Floer homology, and twisted coefficients}

\subsection{Twisted Coefficients} We don't prove the result for all twisted coefficients. Rather, we commit to a specific setting. Namely, we'll suppose that we add a special component to $\Gamma$, denoted $\Xi$, which is assumed to be both null-homologous in $Y$ and in each of the $Y(I)$. For
instance, $\Xi$ might bound a surface which intersects each of the components of $\Gamma$ (if we orient them) algebraically $0$ times. We assume we have a diagram subordinate to $\Gamma \cup \Xi$, and such that in the torus surrounding $\Xi$, the framing is the $0$ framing for the prescribed surface, and the $\be$ curves, and all the $\eta(I)$ curves, are small Hamiltonian isotopes of this framing. We will consider the spectral sequence induced by $\LL$ on the result of $0$-surgery on $Y_{\Xi}$, the result of $0$ surgery on the new curve. We will twist only by the homology class introduced by this surgery. Returning to the Heegaard diagrams, we assume that there is an $\al$-curve and a point on that curve, so that all the isotopes of the $0$-framing curve intersect $\al$ twice near the point, and lie on one side of the $\al$ curve near the point. We call this point, $z'$, and the goal is to arrange it to lie in the boundary of any doubly periodic domain representing the capped surface, and to lie in the same component as $w$ when we have diagrams only containing sets of $\eta$ curves. The prescription above should do this, as the $\al$ curve is necessary to the periodic domain, and the diagrams with all $\eta$ curves decompose into genus $1$ summands with $z'$ outside an small annular neighborhood of the longitude for $\Xi$. We can then define chain maps: 
$$
\widehat{f}_{\mathbb{F}[T,T^{-1}]}({\bf x} \otimes T^{p}) = 
\sum_{{\bf y} \in \T_{\nu^{0}} \cap \T_{\nu^{m}}} \sum_{\{\phi\,|\,\mu(\phi) = 0,\,n_{w}(\phi)=0\}} \#\mathcal{M}(\phi)\big({\bf y} \otimes T^{p- n_{z'}(\phi)}\big)
$$
These maps occur between complexes with twisted coefficients in $\mathbb{F}[T,T^{-1}]$ where the twisting is induced by mapping $H_{2} \ra \Z\cdot F$ where $F$ is the homology class of the surface found by capping the prescribed surface bounding $\Xi$. We can then construct a differential $D^{+}$ as before by summing over sets of immediate successor sequences. As this is just a different additive way of accounting the moduli spaces, we need only check that we have not inadvertently disrupted the cancellations necessary for the theorem. However, these occur in diagrams only including attaching curves from the $\eta(I)$'s. In this case $n_{z'}(\phi) = n_{w}(\phi)$ and it is so there is no change from the $HF^{+}$ case. Thus, we have a spectral sequence with twisted coefficients. This is probably the simplest such version, but it is all we will need.\\

\subsection{Knot Floer Homology} The arrangement of the bouquet allows us to include a component $K$, which will be framed by its meridian, and
such that in this meridian there is a point joined to $w$ by an arc which does not cross any $\al$'s or $\be$'s. As a result, we may put a second point, $z$, on the opposite side of the meridian, and use $w$ and $z$ to encode the knot Floer homology of $K$. If, in addition, this component is null-homologous, and the other components of $\Gamma$ link $K$ algebraically $0$ times, then there is an embedded Seifert surface with boundary $K$ which exists in all the three manifolds $Y(I)$. The knot Floer homology will then exist for each of the pairs $(Y(I), K)$, and will have filtration
in $\Z$. Indeed, the filtration for each $Y(I)$ will be coordinated by this surface, $F$, so that each generator of the knot Floer homology of $(Y(I), K$) will receive the filtration given by 
$$
\mathcal{F}({\bf x}) = \frac{\langle c_{1}\big(\Sps{s}{w}({\bf x})\big), \widehat{F} \rangle}{2}
$$
where the $Spin^{c}$ struture is that induced on $Y_{K}(I)$, the result of Seifert framed $0$-surgery on $K$, and $\widehat{F}$ is the found by capping off the prescribed Seifert surface. As we have been doing heretofore, though not explicitly, we consider the knot Floer groups after 
direct summing over all the $Spin^{c}$ structures on $Y(I)$. \\
\ \\
\noindent The maps in the definition of $D$ will be adjusted as in \cite{Knot} to respect the filtration induced by $K$. Namely, we now
use
$$
\widehat{f}_{\nu^{0}, \ldots, \nu^{m}}( [{\bf x}, 0, j] \otimes \Theta_{0,1} \otimes \cdots \otimes \Theta_{k-1,k}) = 
\sum_{{\bf y} \in \T_{\nu^{0}} \cap \T_{\nu^{m}}} \sum_{\{\phi\,|\,\mu(\phi) = 0\}} \#\mathcal{M}(\phi)[{\bf y}, 0, j - n_{z}(\phi)]
$$
where we assume that $j = \mathcal{F}({\bf x})$. Of course, we should only use those $\phi$ with $n_{w}(\phi) = 0$. Once again, the argument in the preceding sections occurs more or less unchanged, since it depends only on the compactness and formal dimension of the moduli spaces. There are three additional points to check. First, the cancellation lemma when $k \leq 2$. Note, however, that if all the attaching circles are of the $\eta(I)$-variety, then $w$ and $z$ are in the same
component of the complement of these circles in $\Si$. Thus the cancellation occurs as it does for the $HF^{+}$ theory with no change. \\
\ \\
\noindent Second, we need to check that the image of the maps above lies in the subgroup used to define the knot Floer homology. This
will follow from the fact that ${\bf x}$ and ${\bf y}$ as above will be joined by such a homotopy class if an only if:
$$
\mathcal{F}_{I^{k}}({\bf y}) - \mathcal{F}_{I^{0}}({\bf y}) = n_{z}(\phi) - n_{w}(\phi)
$$
This relation ensures that the image is in the subgroup called $CFK^{\infty}$ in \cite{Knot}. It also shows that when $n_{w} = 0$, the maps preserve the filtration, $\mathcal{F}$. Thirdly, we need to check that the various chain homotopies are also filtered. We will give a general argument for why an identity such as the one also hold for the maps in the chain homotopies. Once this identity is established the proof proceeds as before, except for a spectral sequence of filtered complexes. The homological algebra in the appendix generalizes the key lemma to this setting and yields:

\begin{theorem}
Let $\LL$ be a framed link in $(Y, K)$ such that $\mathrm{lk}(L_{s}, K) =0$ for all $s$. For each integer $k$, and surface $F$ spanning $K$ and disjoint from $L$, there is a spectral sequence such that
\ben
\item The $E^{1}$ page is $\oplus_{I \in \{a,b\}^{n}} \widehat{HFK}(Y(I),K, k)$
\item The $d_{1}$ differential is obtained by adding all $\widehat{F}_{I<I'}$ where $I'$ is an immediate successor of $I$ 
\item All the higher differentials respect the dictionary ordering of $\{a,b\}^{n}$, and
\item The spectral sequence eventually collapses to a group isomorphic to \\$\widehat{HFK}(Y(\{cc\ldots c\}, K, k)$. 
\een
\end{theorem}

\noindent In fact, more can be said, since the proposition in the appendix also implies an $r$-quasi-isomorphism on the spectral sequence for
the mapping cone versus the complex $A_{3}$:

\begin{lemma}\label{lem:spec}
For each $r \geq 1$, the $E^{r}$ page of the spectral sequence for $\widehat{HF}(Y)$ computed from 
$\oplus_{k \in \Z} \widehat{HFK}(Y, K, k)$ and arising from the filtration from the knot Floer homology, is quasi-isomorphic 
to the $E^{r}$ page of the complex $X(\mathcal{C})$ using the coherent filtration induced from $K$ on each summand.
\end{lemma}

\noindent {\bf Proof:} The construction of the spectral sequence shows that $X(\{a,b\}^{n-1}\times\{c\})$ is $1$-quasi-isomorphic to
the mapping cone of $X(\{a,b\}^{n-1}\times\{a\}) \lra X(\{a,b\}^{n-1}\times\{b\})$, hence they are $r$-quasi-isomorphic for the
filtration induced by the knot filtration. The mapping cone is the complex $X(\{a,b\}^{n})$, while $X(\{a,b\}^{n-1}\times\{c\})$ can itself
be thought of as the mapping cone of $X(\{a,b\}^{n-2}\times\{ac\}) \ra X(\{a,b\}^{n-2}\times\{bc\})$. This in turn is quasi-isomorphic to
$X(\{a,b\}^{n-2}\times\{cc\})$, and hence $r$-quasi-isomorphic. Continuin down this ladder, we ultimatley arrive at an $r$-quasi-isomorphism with 
$\widehat{CF}(Y(\{cc\ldots c\})$ with the knot filtration. This proves that $X(\{a,b\}^{n})$ is $r$-quasi-isomorphic to $\widehat{CF}(Y(\{cc\ldots c\})$. $\Diamond$\\
\ \\
\noindent To verify the identity above, note that in the corresponding cobordism there is a surface consisting of two copies of the spanning surface, one in each end $\al\eta(I^{0})$ and $\al\eta(I^{k})$, and a neck $K \times I$ joining the boundaries. Since the surface exists throughout the surgery process, this is the boundary of $F \times I$, and is thus null-homologous. Furthermore, the surface in $\al\eta(I^{0})$ will receive the reverse orientation, thought of as the boundary of $F \times I$. Now, $\phi$ determines a $Spin^{c}$ structure on the cobordism, and this $Spin^{c}$ structure
must pair to be $0$ with this surface. Furthermore, $K \times I$ will intersect the subspace determined by $w$ and $\phi$ in $\Si \times \triangle$ transversely and algebraically $n_{w}(\phi) - n_{z}(\phi)$ times. Given the construction of a $Spin^{c}$ structure from this subspace, we see that the pairing of the $Spin^{c}$ structure on the cobordism with the surface $F \cup K \times I \cup - F$ is both zero and 
$$
\langle c_{1}\big(\Sps{s}{w}({\bf y})\big), \widehat{F} \rangle - \langle c_{1}\big(\Sps{s}{w}({\bf x})\big), \widehat{F} \rangle + 2\big(n_{w}(\phi) - n_{z}(\phi)\big)
$$
from which the identity then follows. At each intersection point, the almost complex structure on $\Si \times \triangle$ is reversed, introducing
a change of $2$ in the pairing. Note that this argument does not assume that we have an immediate successor sequence, and thus also applies to the chain homotopies used to construct the spectral sequence from the key lemma. \\ 
\ \\
\noindent We note that the long exact sequence for knot Floer homology from \cite{Knot} follows from this spectral sequence by using a framed link consisting of only one knot. \\
\ \\
\noindent {\bf Notation:} We will mainly be interested in the spectral sequence on $X(\mathcal{C})$ generated by the filtration found from the number of $b$'s in the code for each generator. We will denote the various versions of the complexes discussed above by $\widehat{X}(Y; \LL)$ for the $\hat{\ }$-complex; $X(Y; \LL, M)$ for the twisted coefficient theory, where $M$ is an appropriate module for $\F[H^{1}(Y;\Z)]$; $X(Y, K; \LL, j)$ for the level $j$ knot Floer homology complex and $X^{+}(Y,\LL)$ for the $+$-theory. We will use the same notation for the pages of the resulting spectral sequence, replacing the $X$ by the standard $E$. Thus $E^{k}(X; \LL)$ is the $k^{th}$ page of the spectral sequence for $X(Y;\LL)$. 

\section{Invariance}

\noindent We now consider the complex $X(\mathcal{C})$. First, we show that the particular Heegaard diagram subordinate to the bouquet does not measurably affect the complex. Then we show that the complex obtained from a different bouquet for the same framed link will give us an isomorphic homology. To this end, note that we can build the Heegaard triple underlying $X(\mathcal{C})$ so that the two framings for each component of $\LL$
intersect in only one point, and otherwise do not intersect any other $\be$-curve, \cite{Smoo}. By an argument in \cite{Smoo} proving Lemma 4.5, we 
know that
\begin{lemma}
Any two Heegaard diagrams subordiante to a bouquet $\Gamma$ for $\LL$ can be connected by a sequence of moves chosen from the following:
\ben
\item Handleslides and isotopies among $\seta$.
\item Handleslides and isotopies among $\{\be_{n+1}, \ldots, \be_{g}\}$
\item Isotopies of the framing curves which do not alter the intersection hypotheses
\item Handleslides of the framing curves across $\{\be_{n+1}, \ldots, \be_{g}\}$ 
\item Stabilization
\een
\end{lemma}
\noindent All such moves are presumed to take place in the complement of $w$, and any other marking data, which we have assume to occur in a contractible
subset of the Heegaard diagram. Of these moves, only the third requires elaboration. We surger out the $\be_{i}$ for $i \geq n$. This leaves a diagram in a genus $n$ surface which further decomposes into genus $1$ sub-diagrams. We view this as $n$ handles attached to a sphere, and locate $w$ in the sphere. Now the framings determine a basis for $H_{1}$ of each genus $1$ component. As such we can take the $\Z$ covering along one of the circles to obtain a cylinder with a single puncture in each translate, corresponding to the rest of the diagram. Any other isotope of the other framing curve, which intersects once, lifts to an infinite line. Therefore, we can isotope any two such lifts, one to the other in each translate, without increasing the intersection number, and without crossing the puncture. Note that since all of these are considered to be $\be$-curves we cannot isotope a non-framing $\be$-curve to intersect a framing one. Since the framings are a spine for each punctured genus $1$ component, we may assume that any move involving only non-framing $\be$'s occurs away from a neighborhood of the framing curves. We need these observations to ensure that $X(\mathcal{C})$ is still a chain complex after handleslides and isotopies.\\
\ \\
\noindent In addition to checking these moves, there are two other invariances to confirm. The first is to see that the particular bouquet chosen
does not affect the quasi-isomorphism type of $X(\mathcal{C})$. We will only verify that there is a $1$-quasi-isomorphism. We do this at the end of the section. Second, since each $\be$-curve, and framing curve, needs to be duplicated, we had to introduce points $1_{i}$ and $2_{i}$ on each of these curves. We will verify that we can slide these points past the intersections of $\be_{i}$ (or the framing curves) with other curves. For the non-framing $\be$-curves this merely means intersections with the $\al$-curves. For the framing curves this will mean both $\al$-curves and the other framing curve. We do this when we verify the isotopy invariance.\\

\subsection{Relabelling} We start with the simplest invariance. $X(\mathcal{C})$ is isomorphic to the chain complex $X(\mathcal{C}')$ obtained
from the same link with the same framings, but with the order of the link components changed. This follows from the symmetry in the definition
of the complex.\\
 
\subsection{Change of $J$'s} In addition to the Heegaard data we must also prove that the choice of almost complex structure data on $Sym^{g}(\Si)$
does not affect the complex materially. So far, in these notes, we have ignored all such questions, refering instead to \cite{Seid} or \cite{Fooo}.
Here we will onnly give a brief synopsis of the argument. The almost complex data for our diagrams needs to be a map $J: \mathcal{M}(\triangle) \times \triangle \ra \mathcal{U}$. Here, $\mathcal{M}(\triangle)$ is the moduli space of conformal structures on a polygon with $3$ or more sides. We choose
this to be the largest polygon appearing in the complex, which is only dependent on the number of components of $\LL$. $\mathcal{U}$ is the set of nearly symmetric almost complex structures on $Sym^{g}(\Si)$ used in \cite{Hom3}. Furthermore, we must add conditions in neighborhoods of each intersection point for the all the totally real tori involved. Namely we must specify that the almost complex structures limit to specified sets
in these neighborhoods so that we can glue moduli spaces near these points. We will not be specific about this. \\
\ \\
\noindent We now consider a (generic) path $J^{(\tau)}$ of such families of almost complex structures for $\tau \in [0,1]$. We define a map
$$
\mathcal{J}(\xi) = \sum_{J \in \mathcal{S}}\  \sum_{I = \codeS{0}{k} = J} \mathcal{J}_{\codeS{0}{k}}(\xi)
$$
where the sum is over immediate successor sequences, and $\mathcal{J}_{\codeS{0}{k}}$ is defined by
$$
\mathcal{J}_{\codeS{0}{k}}(\xi) = \sum_{\zeta} \sum_{\psi} \#\big(\bigcup_{\tau\in[0,1]} \mathcal{M}(\psi, J^{\tau})\big)\zeta
$$
where $\mu(\psi) = -1$. Note that $\xi$ and $\zeta$ may come with more information: indices, twisted coefficients, etc. In that case
we modify the map as we did with the maps in the spectral sequence. Note that the count is the number of points in a union of moculi spaces, one for each $\tau \in [0,1]$. Since this family depends on a $1$-dimensional paramter, a $\mu = -1$ homotopy class will have $0$ formal dimension for its moduli space. We note]that when $k =1$, we obtain the map for changing the path of almost complex structures on the Heegaard-Floer homology. \\
\ \\
\noindent We now consider a map built in exactly the same manner, but with $\mu = 0$ homotopy classes. These have moduli spaces with compactifications, and the only contributions to the sum are those from the broken boundary components and when $\tau =0$ or $\tau = 1$. If the underlying homotopy class is that of a $k+1$-gon, then there are two types of degeneration: 1) when the division occurs uses the $\al$-edge, and 2) when the division uses $\eta(I)$ edges only. In the first case we obtain $D_{1} \circ \mathcal{J} + \mathcal{J} \circ D_{0}$, i.e. all the variation in the almost complex structure occurs in one of the polygons. This includes the case when we divide along $\al$ and $\eta(I^{0})$ or $\eta(I^{k})$ and $\al$, where one of the maps is the chain map for the Heegaard-Floer homology of the respective end. On the other hand, the cancellation lemma
still applies to the second case, since the Maslov index was found by splicing triangles, and this will not change as we alter $J$. Thus, $\mathcal{J}$
will be a chain map. Furthermore, we can reverse the path to obtain another chain map. Since both of these induce the standard chain maps, for varying the almost complex structure, on $CF(Y(I))$, they induce isomorphisms on the $E^{1}$-pages of the spectral sequences which are inverses of each other.\\
\ \\
\noindent This argument is not strictly correct, but contains the elements necessary for improvement. The main problem stems from not being explicit
about the almost complex data and the compactifications, but I am not a trustworthy guide to this arena. In addition, there are difficulties with when to allow the paths of almost complex structures near the intersection points to vary. Nevertheless, as in \cite{Hom3}, \cite{3Man}, these concerns should be assuageable, so the interested reader should now look at \cite{Seid}, for example.\\

\subsection{Under Handleslides}

\subsubsection{Handlesliding the $\al$-curves} Let $\al'$ be the result of a single handleslide among the $\al$-curves of the diagram. This should not interrupt any assumptions about the placement of marked points and their adjacency. We consider a map formed
$$
B_{\mathcal{S}}(\xi) = \sum_{J \in \mathcal{S}}\  \sum_{I = \codeS{0}{k} = J} B_{\codeS{0}{k}}(\xi)
$$
where the sequence is an immediate succesor sequence and $B_{\codeS{0}{k}}$ counts holomorphic $k+3$-gons with boundary on $\al, \eta(I^{0}), \ldots, \eta(I^{k}), \al'$ thought of as a map between $CF(\al, \eta(I^{0}))$ and $CF(\al', \eta(I^{k}))$. The rest of the intersection points, including $\al\al'$ consists of the canonical generator $\Theta^{+}$ as before. Note that we now have generators on the right and left. As usual we consider one dimensional moduli spaces of such things, and then look at the boundary of their compactifications. If the division occurs in the $\al$ edge, and some other edge then 1) if the other edge is $\al$' the moduli space cancels with some other space due to the closure of $\Theta^{+}_{\al\al'}$, or 2) the other edge is $\eta(I^{j})$ for some $j$. In the latter case we obtain a portion of the map $B \circ D_{\al}$ for the differential above. If the division occurs including $\al$' and an $\eta(I^{l})$ then we obtain $D_{\al'} \circ B$. The other divisions occur between two sets of attaching circles, $\eta(I^{i})$ and $\eta(I^{j})$. One of the polygons resulting from this divison has the form in the cancellation lemma. Adding over all the immediate successor sequences shows that the contribution of the polygons which can fulfill the same role will be zero since we have not alterred any of the $\eta(I^{j})$ curves. Thus $B \circ D_{\al} + D_{\al'}\circ B = 0$ and $B$ is a chain map. The map induced on $E^{1}$ is the standard handleslide map for the $\al$-curves. As a result, it induces an isomorphism of the $E^{1}$ pages and thus an isomorphism of every page. \\

\subsubsection{Handleslides of the non-framing $\be$-curves} 
\begin{center}
\begin{figure}
\includegraphics[scale=1.0]{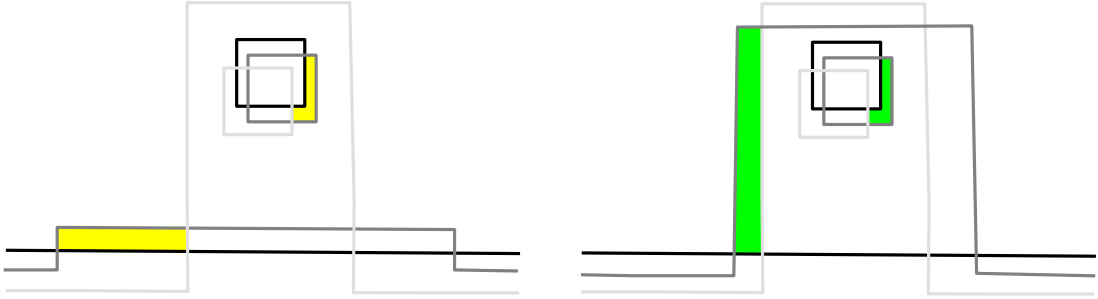} 
\caption{These diagrams depict the triples used in verifying cancellation for the handleslide of a non-framing $\be$ over another non-framing $\be$. On the left, we have the curves for $I^{0}, I^{1}, I^{2}$, from darkest to lightest respectively, for a code sequence $0\ldots;0 \ra 1\ldots;0  \ra 1\ldots;1$ with holomorphic triangles in yellow. The right depicts the diagram for $0\ldots;0 \ra 0\ldots;1 \ra 1\ldots;1$, with holomorphic triangles in green. Note that these triangle have the same vertex for $\Si_{I^{0}I^{2}}$. All other local diagrams is of one of the types in Figures \ref{fig:TypeIII}, \ref{fig:TypeII}}\label{fig:Handleslide1}
\end{figure}
\end{center}

We consider the case of sliding a $\be$ not involved in a framing over another $\be$
curve not involved in the framings. For the complex $X(\mathcal{S})$ we duplicate each $\be$-curve multiple times. However, this is canonically specified in the annular neighborhood of the $\be$-curve. As we noted above, isotoping in this annular neighborhood, but not changing the intersection data, corresponds to changing the underlying complex structure on $\Si$, which in turn is used to characterize admissible $J$'s. This argues
that we do not need to consider handlesliding each individual duplicate; we will instead do this all at once. The pattern of the proof is similar
to the previous case, so we will only identify the differences. Our codes now end in $\ldots0;1$, for example, where a $1$ after the semi-colon
indicates that we use $\be_{i}'$ to form the Hamiltonian isotopes, where $\be_{i}'$ is the result of handlesliding $\be_{i}$ over $\be_{j}$ ($i,j > n$). A $0$ after the semi-colon indicates that we use $\be_{i}$ as usual. \\
\ \\
\noindent The maps now count holomorphic $k+2$-gons with boundary on $\al, \eta(I^{0}), \ldots, \eta(I^{k})$, $\mu = 0$ and mapping through
the standard generators for $\eta(I^{i})\eta(I^{i+1})$. Furthermore, there must be a pair, $I^{i}$ and $I^{i+1}$, where the change in code is 
$\ldots;0 \ra \ldots;1$. To prove that this is a chain map we consider the $1$-dimensional moduli spaces. The degenerations occur as usual, and the only question is whether we still have a cancellation lemma. That $\mu = k - 2$ for $\eta(I^{0}), \ldots, \eta(I^{k})$ for $k > 2$ follows
as before, noting that there is a small holomorphic triangle in the picture for $\eta(I^{0}), \eta(I^{i}), \eta(I^{i+1})$ where the $0$ after the 
semi-colon changes to $1$. The doubly periodic domains are now larger, but still do not affect the dimension. Furthermore, $[S]$ cannot be added or subtracted without changing $n_{w}$ and the dimension of the moduli space. Thus, sequences where we change the element after the semi-colon will have a cancellation lemma, as do sequences where we do not (by the usual cancellation lemma). This leaves the $k=2$ case. For sequences without the handleslide, the usual lemma applies. When we do have the handleslide, we have a code change of the form $0 \ldots;0$, $0 \ldots;1$, $1 \ldots;1$, or $0 \ldots;0$, $1 \ldots;0$, $1 \ldots;1$. Checking the two genus $1$ components involved shows that the triangles for these changes will have the same image, and thus cancel, see Figure \ref{fig:Handleslide1}. Thus the usual argument establishes that the map is a chain map.
\noindent We have built the map to induce on the $E^{1}$ pages the maps from \cite{Hom3} used to show the handleslide invariance of the Heegaard-Floer
homology. Since these maps occur between Heegaard-Floer homologies of the $Y(I)$, this shows that they induce isomorphisms on every page of the
spectral sequence. One imagines that we can do this in the category of filtered chain maps, up to filtered chain homotopy -- and, indeed, we will prove
a composition property below which will also apply here -- but we do not do not need this. \\

\subsubsection{Handleslides of framing curves over non-framing curves} 
\begin{center}
\begin{figure}
\includegraphics[scale=1.0]{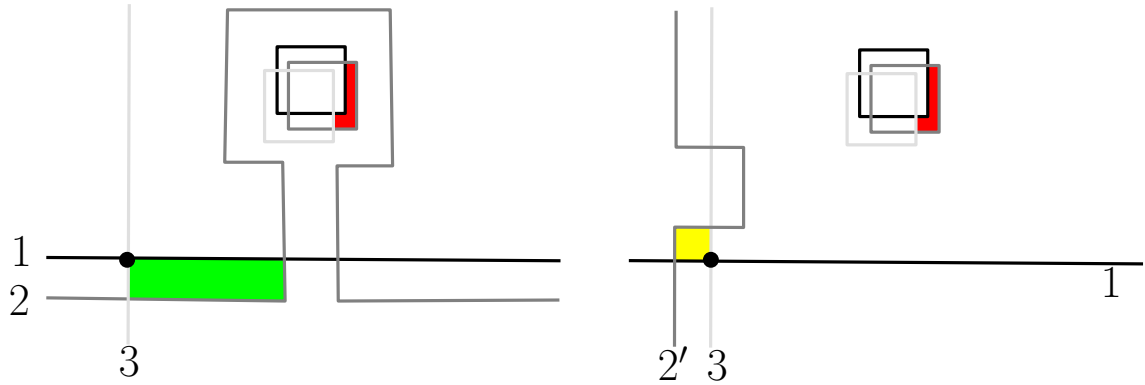} 
\caption{These diagrams depict the triples used in verifying cancellation for the handleslide of a framing $\be$ over a non-framing $\be$. In each we have the curves for $I^{0}, I^{1}, I^{2}$, from darkest to lightest respectively, with the $b$ framing in the light gray intersecting the other two curves in one point. On the left we have the diagram for a code sequence $0\ldots;0 \ra 1\ldots;0  \ra 1\ldots;1$ with holomorphic triangles colored. The right depicts the diagram for $0\ldots;0 \ra 0\ldots;1 \ra 1\ldots;1$. Note that these triangle have the same vertex for $\Si_{I^{0}I^{2}}$. All other local diagrams is of one of the types in Figures \ref{fig:TypeIII}, \ref{fig:TypeII}}\label{fig:Handleslide2}
\end{figure}
\end{center}

The argument is the same as all the arguments we have been presenting, so we will simply say how to adjust the codes and how to obtain the cancellation. The code $\ldots;0$ will mean use the curves before the handleslide, while
$\ldots;1$ will mean use those after the handleslide. There are two possibilities depending on which framing curve we handleslide. Here we will assume
that it is the $b$-curve; the argument for the $a$-curve is boringly similar. Cancellation follows as in the previous case, except we need to
examine the special triangles that arise at the end. There are two sets of local pictures for $0 \ldots;0$, $0 \ldots;1$, $1 \ldots;1$, or $0 \ldots;0$, $1 \ldots;0$, $1 \ldots;1$ depending on whether the first index corresponds to the framing curve being slid, or not. \\
\ \\
\noindent If it doesn't, then either the index which does is $1$, and the argument is, mutatis mutandis, that of the previous subsection, or it is $0$, and then the argument is as in the standard cancellation lemma, since this corresponds to using the $a$-framing curve, which has not been alterred. If it does, then in one case we alter the framing before the handleslide, and the handleslide has no effect. In the other, we slide then alter the framing. A local calculation shows that these still have the same image, and therefore cancel, see Figure \ref{fig:Handleslide2}
\noindent We conclude that the underlying map is a chain map, and that it induces the handleslide isomorphisms on the $E^{1}$-page. For those $Y(I)$
using the $a$-curve, the $E^{1}$ map is just the identity. It is invertible, using the map from the reverse handleslide. \\

\subsection{Under Isotopies}

\subsubsection{A useful observation} The complexes $X(\mathcal{S})$ are built out of Heegaard multi-tuple maps. It is therfore convenient to
change the maps induced on the Heegaard-Floer homlogies by isotopies into maps built out of Heegaard triples. This is not obviously possible at the
chain level, but it possible at the level of homology. For example, suppose $(\Si, \al, \be)$ is a Heegaard diagram for $Y$, and we perform
a small Hamiltonian isotopy of $\be_{i}$ to obtain $\be_{i}'$. Let $\ga_{i} = \be_{i}'$ and let $\ga_{j}$ be a small Hamiltonian isotope
of $\be_{j}$. Of course, the idea is that the latter are somehow ``smaller'' isotopes, namely occurring in prescribed annular neighborhood, whereas the former can be quite larege. Also, the latter will only be allowed to intersect with their model twice, whereas no such requirement is made on the former. Nevertheless, we have (where $HF$ stands in for any flavor of Heegaard-Floer homology):

\begin{lemma}
The map on $HF(Y)$ by the isotopy of $\be_{i}$ to $\be_{i}'$ can be written as a composition of Heegaard triple maps.
\end{lemma}

\noindent{\bf Proof:} Consider the isotopy taking $\ga_{i}$ back to $\be_{i}$ (really a small Hamiltonian isotope intersection $\be_{i}$ twice only). 
According to section 8 of \cite{Hom3}, the chain map $\Gamma^{-1}$ induced by this isotopy on $(\Si, \al, \ga)$ has the property that
$$
\Gamma^{-1} \circ f_{\al\be\ga} + f_{\al\be\ga'} \circ \Gamma_{\be\ga\ga'} = \partial \circ H + H\circ \partial
$$
where $\ga'$ is the result of isotoping the $\ga$-curves, i.e. $\ga'$ is just small isotopes of $\be$. Thus $f_{\al\be\ga'}$ will induce
an isomorphism on the chain complexes whose top level in the area filtration is the identity. Meanwhile $\Gamma_{\be\ga\ga'}$ wil take
$\Theta^{+}_{\be\ga}$ to $\Theta^{+}_{\be\ga'}$ on homology since these are the only generators in their grading. Thus, on the homology
we obtain that $\Gamma^{-1} \circ F_{\al\be\ga} = F_{\al\be\ga'}$. Since the last map is an isomorphism, both of the $F$-maps are isomorphisms and we can write $\Gamma = F_{\al\be\ga} \circ F^{-1}_{\al\be\ga'}$. $\Diamond$ \\
\ \\
\noindent We can then construct a map on $X(\mathcal{S})$ which induced on $E^{1}$ the map $F_{\al\be\ga'}$, where we make a small isotope to
each $\be$ curve. This is done by specifying a code $\ldots;0$ or $\ldots;1$ with an additional entry. The $0$ means use the $\be$ and the framing curves as models for the small isotopes, the $1$ says use the $\ga'$ as models. We then define
$$
B_{\mathcal{S}}(\xi) = \sum_{J \in \mathcal{S} \times\{1\}}\  \sum_{I = \codeS{0}{k} = J} B_{\codeS{0}{k}}(\xi)
$$
where $B_{\codeS{0}{k}}$ counts holomorphic polygons as in the differentials defined above. Drawing local pictures shows that we still have the
cancellation lemma in this case. So the standard argument ensures that this is a chain map inducing the maps $F_{\al\be\ga'}$ on the $E^{1}$
page, and thus inducing isomorphisms on every page.\\
\ \\
\noindent We now show how to obtain a similar picture for $F_{\al\be\ga}$. The composite of this with the inverse of the preceding shows that
isotopies also induce isomorphisms on every page. The argument for isotoping $\al$-curves is similar; indeed, it has fewer cases and is thus simpler.

\subsubsection{Isotopies within an annular neighborhood} 
\begin{center}
\begin{figure}
\includegraphics[scale=1.0]{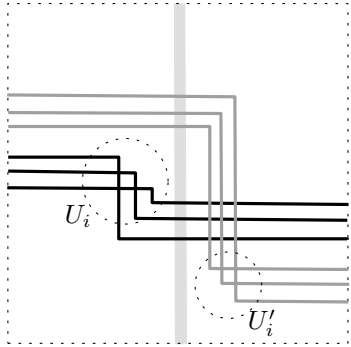} 
\caption{Pushing the intersection region $U_{i}$ past an $\al$ curve. The dark curves correpsond to using $\ldots;0$, whereas the lighter curves iccur for $\ldots;1$. Note that without the $\al$ curve, and including the remainder of the annulur neighborhood of the $\be$-curve, we obtain a picture similar to the right side of Figure \ref{fig:splice2}.}\label{fig:SpecPtIso}
\end{figure}
\end{center}

\begin{center}
\begin{figure}
\includegraphics[scale=1.0]{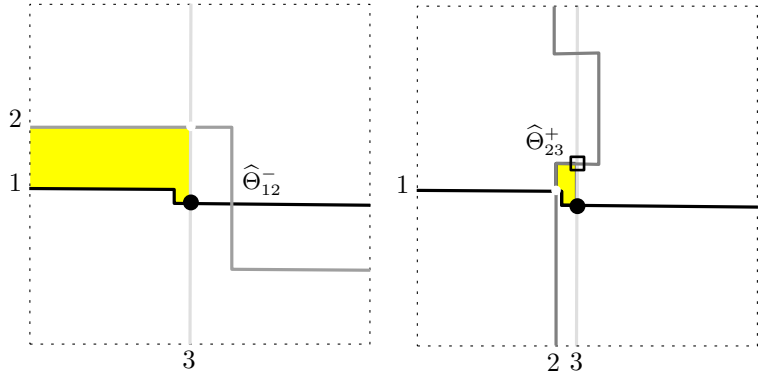} 
\caption{On the left we have the code sequence $0\ldots;0 \ra 0\ldots;1 \ra 1\ldots;1$ corresponding to two small isotopes of $F_{a}^{i}$ and a single, light gray copy of $F_{b}^{i}$, which we have slide the past the region $U_{i}$. There is a holomorphic triangle, although the intersection point $\Theta^{+}$ occurs outside the diagram, along the curves. On the right, corresponding to the code change $0\ldots;0 \ra 1\ldots;0 \ra 1\ldots;1$ there are two small isotopes of $F^{i}_{b}$ and one of $F_{a}^{i}$. Now the holomorphic triangle is in the local picture, but both holomorphic triangles have the same intersection point.}\label{fig:SpecPtTrian}
\end{figure}
\end{center}

We will first check that we can move the points $1_{i}$ and $2_{i}$. Since these are symmetric, we will only verify the case for $1_{i}$ when $\be_{i}$ is non-framing. When $\be_{i}$ is framing, the same argument will apply if $1^{a}_{i}$ or $1^{b}_{i}$ move past an $\al$-curve. So we will only verify the invariance when we move $1^{a}_{i}$ past the intersection of $F_{i}^{a}$ and $F^{b}_{i}$. When moving $1^{b}_{i}$, the argument is formally similar. \\
\ \\
\noindent We enhance the code space to include another entry: $\ldots;0$ or $\ldots; 1$. We interpret these using: $0$ means use $1_{i}$
on one side of an $\al$ curve, and $1$ means use $1_{i}$ on the other side of the $\al$-curve. Now we construct the map $B$ as above, using
immediate successor sequences with this enchanced entry. To obtain the chain map, we assume that in the code sequence, the last entry must change. When we consider one dimensional moduli spaces we may apply the $A_{\infty}$-relation. Divisions along the $\al$ edge produces two pieces, only one of which can contain the change in the last entry. That piece corresponds to an term in the chain map. The other pieces occurs in the chain complex for the pieces with all $0$ or all $1$ in the last entry. Thus, we need only verify the cancellation lemma to obtain a chain map. Divisions for the polygons in our diagrams that divide off an all $\eta$-portion come in two forms: 1) when the last entry is always $0$ or $1$ or 2) when it switches. In the former case, the usual cancellation lemma applies, as both the diagram before and the diagram after the slide have the right structure. There is no $\al$-curve in the all $\eta$-diagram, after all. So we concentrate on those diagrams where a switch does occur. In this case we can assume that the diagram is as in the local picture Figure \ref{fig:SpecPtIso}. This diagram is the normal diagram for a genus $1$ component near a grouping of intersection points, except that we have partitioned the curves into two sets. However, in the all $\eta$-diagram the partition is irrelevent, so the cancellation lemma still applies. 
\noindent Examining the map on the $E^{1}$ page, we see that we need to count triangle with an $\al$-edge and only the last entry changes. This
is easily seen to be equivalent to taking small Hamiltonian isotopes of all the $\be$-curves for the code $I$, and is thus an isomorphism on $HF(Y(I))$. Since the $E^{1}$-map is a direct sum of these isomorphisms, it is also an isomorphism. Therefore, as this is induced from a filtered chain map, all the higher pages are also isomorphic.\\
\ \\
\noindent A very similar argument applies when sliding $1_{i}^{a}$ past the intersection of $F_{a}^{i}$ with $F_{b}^{i}$. We will only indicate
how cancellation occurs. Note that the $E^{1}$-map will be as above, since it arises from the change in the location of $1_{i}^{a}$ without alterring
the code in any other spot. Once we have cancellation, we will be done. Again, if the last entry in the enhanced code doesn't change, the usual cancellation lemma applies. If it does change, and $k > 2$, we construct a $\mu = k-2$ $k+1$-gon, which, as usual, is the only possible
$k+1$-gon for the map. Since it has the wrong dimension, this will suffice. The diagrams look the same as in the usual argument, except for the one triangle, $\eta(I^{0}), \eta(I^{j}), \eta(I^{j+1})$, where the last entry changes. In this case, either we have three small isotopes of each other, and we can find a small holomorphic triangle, or two of the curves are small isotopes, corresponding to $F_{a}^{i}$, and the third curve is an isotope of $F_{b}^{i}$. Again there is a small holomorphic triangle. Thus we need only consider the $k=2$ case, which again occurs for the codes
$0\ldots;0 \ra 0\ldots;1 \ra 1\ldots;1$ and $0\ldots;0 \ra 1\ldots;0 \ra 1\ldots;1$. If the first coordinate does not correspond to $F_{a}^{i}$, then 
the local diagram for the $i^{th}$ component is just three small Hamiltonian isotopes of either $F_{a}^{i}$ or $F_{b}^{i}$, while that for the first entry that changes is either a two copies of $F_{a}^{j}$ and a copy of $F_{b}^{j}$ or vice versa, see Figure \ref{fig:TypeIII}. In either case, the image of the map is the canonical generator $\widehat{\Theta}_{1\ldots;1}$. If it does correspond, then there are two diagrams to consider, see Figure \ref{fig:SpecPtTrian}. For the first sequence, we have two small isotopes of $F_{a}^{i}$ followed by a copy of $F_{b}^{i}$, and three small isotopes of every other curve, for which the image is $\widehat{\Theta}_{1\ldots;1}$, or a copy of $F_{a}^{i}$ followed by two copies of $F_{b}^{i}$, for which the image is the same. Thus we still have cancellation, and so we are done.

\subsubsection{Isotopies of non-framing $\be$'s} 
\begin{center}
\begin{figure}
\includegraphics[scale=1.0]{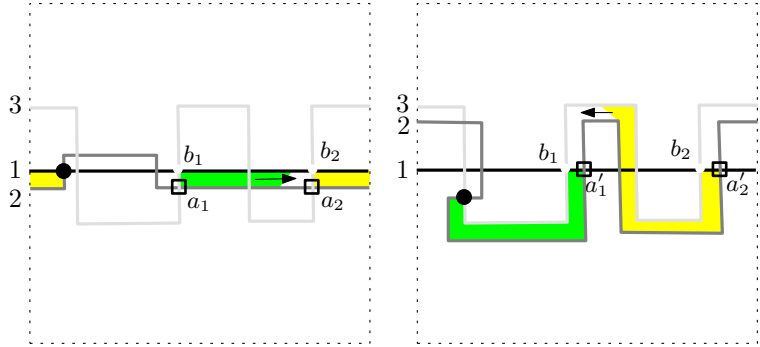} 
\caption{The darkest curve is the $\be$ to be isotoped. The left corresponds to the code sequence $0\ldots;0 \ra 1\ldots;0 \ra 1\ldots1$ for a small Hamiltonian isotopy introducing four intersection points. On the right we have the diagram for this isotopy corresponding to the code sequence $0\ldots;0 \ra 0\ldots;1 \ra 1\ldots1$. As usual the curves go from darkest to lightest as the index increases. On the left, $a_{1} + a_{2}$ is the generator for the $S^{1} \times S^{2}$ factor coming from this pair. It's image is $b_{1} + b_{2}$, the canonical generator for the first and third curve, using the yellow and green holomorphic triangles, the green one extends all the way around $\be_{i}$ until arriving at $\hth{i+1}{i+2}$, indicated by the black disc. On the right, the black disk is still $\hth{i+1}{i+2}$, but now $a_{1}' + a_{2}'$ is the canonical generator for its pair, whose image is $b_{1} + b_{2}$. Thus the two sequences will have the same image, as all other diagrams have the form in Figures \ref{fig:TypeIII}, \ref{fig:TypeII}.}\label{fig:Isotopy1}
\end{figure}
\end{center}

If we choose a $\be_{i}$ with $i > n$, and make a small Hamiltonian isotopy, we arrive at a curve
$\be_{i}'$ which intersects $\be_{i}$ an even number of times. We will consider the map built like $B$ above, except with the following addition to 
the code space: $\ldots; 0$ indicates we use the curve $\be_{i}$ as a model, whereas $\ldots; 1$ indicates we use $\be_{i}'$ as a model. Of course, there is an annular neighborhood of $\be_{i}$ containing $\be_{i}'$; however, we will assume that the annular neighborhoods used when each is a model intersect in the same manner, combinatorially, as the curves themselves. In other words, the small Hamiltonian isotopes arising from the codes occur in much smaller annular neighborhoods than that surrounding $\be_{i}$ and $\be_{i}'$. \\
\ \\
\noindent As in all our arguments, the key issue is whether there is a cancellation for polygons labelled with all $\eta$ attaching circles. If
the codes in the $\eta(I)$'s all end in $0$ or $1$, then the usual cancellation lemma applies. Assume we have a polygon where the final element
changes from a $0$ to a $1$. If $k > 2$, then there are small triangles which can be used to construct a $\mu = k-2$ polygon. In fact, in the triangle
$\eta(I^{0}), \eta(I^{j}), \eta(I^{j+1})$ where the last entry changes, there are distinct small triangles for each of the intersection points
between $\be_{i}$ and $\be_{i}'$. In the next triangle (or previous) triangle, there are also many such small triangles. As with the lemma \ref{lem:workhorse}, all of these will have the same effect on splicing. Thus, we are back to the $k=2$ case, and the code changes $0\ldots;0 \ra 1\ldots;0 \ra 1\ldots;1$ and $0\ldots;0 \ra 1\ldots;0 \ra 1\ldots;1$. We draw the local diagrams to see that the corresponding triangle maps will have the same image, see Figure \ref{fig:Isotopy1} as an example. Of course, we already know this for the genus $1$ component where the framing curve changes, so we only require confirmation in the annular region near $\be_{i}$ and $\be_{i}'$. Obviously, the $E^{1}$-map on each $\HomT{Y(I)}$ will be the triple map $F_{\al\be\ga}$ for the isotopy, described above.
\subsubsection{Isotopies of framing $\be$'s} 
\begin{center}
\begin{figure}
\includegraphics[scale=1.0]{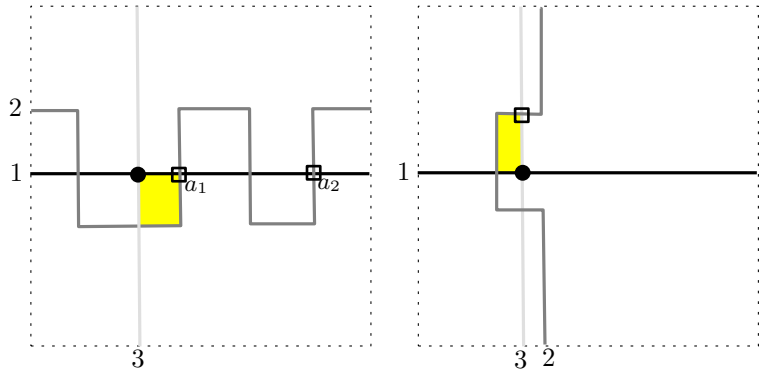} 
\caption{When the isotopy occurs on a framing curve, there are several possibilities. This depicts an isotopy of $F_{a}^{i}$, which either can be isotoped before channging to $F_{b}^{i}$, or we can switch to $F_{b}^{i}$ first, in which case the isotopy plays no role, as on the right. We are already familiar with the diagram on the right, but on the left we see that the canonical generator for the first two curves is $a_{1} + a_{2}$. It's image, through the canonical generator for the second pair, is the same as on the right, since only one of $a_{1}, a_{2}$ abuts a holomorphic triangle. The homotopy class for the other genrator will have negative multiplicities. The other possibilities, e.g. isotoping $F_{b}^{i}$ follow similarly.}\label{fig:Isotopy2}
\end{figure}
\end{center}

Consider the code space with the extra entry with $\ldots;0$ indicating the unisotoped $\be_{i}$ for $i \leq n$, while $\ldots;1$ indicates the isotoped copy, $\be'_{i}$. We note that cancellation follows as before for $k \ge 2$, using the $\mu =0$ small triangles which necessarily arise in the $\hat{\ }$-version of the Heegaard triple map for the isotopy. This reduces us to the $k=2$ case for the code sequences $0\ldots;0 \ra 1\ldots;0 \ra 1\ldots;1$ and $0\ldots;0 \ra 1\ldots;0 \ra 1\ldots;1$. If the other code which changes is not that for $\be_{i}$ then the cancellation result follows from the same considerations as in the isotopies of non-framing $\be$ curves. If the code is the same as the curve which is isotoped, then we use a different set of diagrams, see Figure \ref{fig:Isotopy2}. Nevertheless, the the cancellation result still holds. The $E^{1}$ portion of this map is the $F_{\al\be\ga}$ maps as above, using the appropriate $\ga$: namely, the isotoped $\be$ for the $a$ framing, and the usual $b$-framing. 
\subsection{Under Stabilization}

\noindent Stabilization of the Heegaard diagram subordiante to $\Gamma$ consists of connect summing the Heegaard diagram with a genus $1$ diagram consisting of an $\al$-curve and a $\be$-curve intersecting in a single point and representing $S^{3}$. All the curves in $\eta(I)$ consist of small Hamiltonian isotopes of the $\be$-curve. For a $k+1$-gon consisting of $\eta(I)$ decorations only, this new component is exactly like the other
genus $1$ components coming from $\be$ curves not involved in the framings. By the gluing lemma, we can remove the component without changing the Maslov index or moduli space. Furthermore, since the $\al$-curve intersects each of the $\be$-curves once, we can find a $k+1$-gon for any
sequence $\al, \eta(I^{1}), \ldots, \eta(I^{k})$ which admits a moduli space of dimension $k-2$, found from splicing holomorphic triangles for $\al, \eta(I^{j}), \eta(I^{j+1})$. Once again we can use the gluing lemma. As a result, the maps for any immediate successor sequence for the stabilized
diagram is precisely the same as that for the unstabilized diagram, except that we add the single intersection point from the genus $1$ component for $\al, \eta(I^{j})$ for each intersection point. Thus $X(\mathcal{C}) \cong X'(\mathcal{C})$ under stabilization.

\subsection{Under changes of bouquet} 
\begin{center}
\begin{figure}
\includegraphics[scale=0.8]{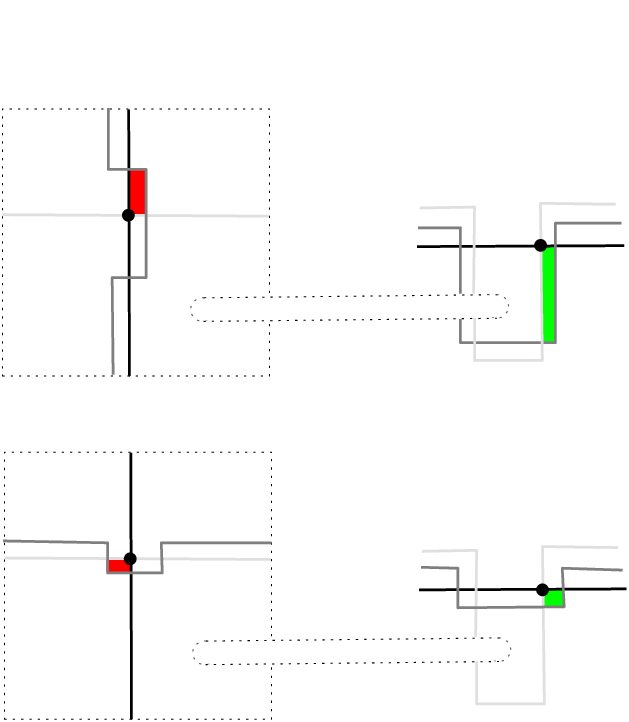} 
\caption{The dotted square is a genus $1$ summand to our Heegaard diagram, including two framing curves intersecting in a point, along with their annular neighborhood and isotopes. It joins the remainder of the diagram by the dotted tube. The triples depict the process of sliding a non-framing $\be$ across the connect sum region for such a genus $1$ component. Done in steps, this would require us to slide twice over whichever framing curve was present, but we will simply jump to the end. The top depicts the diagram for a code change of $0\ldots;0 \ra 0\ldots;1 \ra 1\ldots;1$ where the framing curve which changes is the one in the genus $1$ component. The bottom depicts the code change $0\ldots;0 \ra 1\ldots;0 \ra 1\ldots;1$ for the same entries. As usual, these code changes correspond to triangles with the same image in $\Si_{13}$, given by the colored homotopry classes of triangles.}\label{fig:bouquet}
\end{figure}
\end{center}

In \cite{Smoo}, lemma 4.8, it is shown that given another bouquet, $\Gamma'$, for $\LL$, there
are two diagrams, one subordinate to $\Gamma$ and the other subordinate to $\Gamma'$ such that one can be obtained from the other
by handleslides of attaching circles. Indeed, this only requires handleslides of $\be$'s, but will include a handleslide across a framing curve.
Since we have already dealt with isotopies and handleslides over non-framing curves let us consider the result of handlesliding past the framing curve. In fact, in the proof of lemma 4.8 we slide twice over the $a$-framing curve, so it suffices to consider the following situation: $\be_{i}$ starts on one side of the connect sum neck for a genus $1$ summand containing an $a$ and $b$ framing curve, while $\ga_{i}$ is the same curve on the other side of the connect sum neck. We attach an extra entry to the codes as before: $\ldots; 0$ means use $\be_{i}$ as a model for small Hamiltonian isotopes, while $\ldots;1$ means use $\ga_{i}$. Once again the crucial question is whether we can attain cancellation. For immediate successor sequences where all the codes end in $0$ or $1$, cancellation follows from the original cancellation lemma. For sequences where the last entry changes from $0$ to $1$, we follow the pattern as above using the local picture in Figure \ref{fig:bouquet} to find a small holomorphic triangle to splice for the change from $\be_{i}$ to $\ga_{i}$. Note that the diagram in the genus $1$ summand is one of the standard diagrams which can be removed without changing the Maslov index. The remaining diagram looks as if we have just isotoped $\be_{i}$ and so the result will hold. Note that even though we have a large doubly periodic domain, adding or subtracting it will still introduce negative multiplicities to the spliced $k+1$-gon. This reduces us to the $k=2$ setting and consideration of code sequences $0\ldots;0 \ra 1\ldots;0 \ra 1\ldots;1$ and $0\ldots;0 \ra 0\ldots;1 \ra 1\ldots;1$. We must show that these have the same image. This follows again from the local pictures for the curves $\be_{i}$ and $\ga_{i}$ and the genus $1$ summand. Therefore cancellation occurs. 

\noindent The $E^{1}$ map induced by this chain map is an isomorphism since we can reverse the process and obtain another map for which composition
will be an isomorphism at the $E^{1}$ level. Thus all the $E^{r}$ pages are isomorphic. Note that for some of the theories, namely the twisted coefficients theory and the knot Floer homology, it is possible that handleslides of this type might pass over the marking data. We show now that this need not be the case. Basically, we surger out all the $\be$ curves except those representing $\be_{i}$ and the framing on the knot $K$ used in the two theories. This leaves a genus $2$ diagram, with $\be_{i}$ the meridian of one of the components. Now $\ga_{i}$ is a meridian of the same component, since it is not allowed to intersect either of the framings on $K$, which span $H_{1}$ of the torus summand. $\be_{i}$ can be slid to $\ga_{i}$ around a path away from the genus $1$ component, unless $w$ is placed in the way. However, our assumption on $w$ is that there is an unobstructed path to one or the other of the pair of framings. Thus $w$ must effectively lie in the genus $1$ component and cannot get in the way of
moving $\be_{i}$ to $\ga_{i}$. Since this movement is unimpeded the results from the previous paragraphs apply, and the complex $X(\mathcal{C})$
is $1$-quasi-isomorphic to the complex $X'(\mathcal{C})$.  

\section{Maps on the spectral sequences}

\subsection{Maps from one handle attachments} Let $W$ be the four dimensional cobordism obtained by taking $Y \times I$ and attaching a one handle
to $Y \times \{1\}$. Let $Y' \cong Y \# S^{1} \times S^{2}$ be the resulting boundary. Let $\LL \subset Y$ be an $n$-component link equipped with triads of framings. The addition of the one handle does not affect the framed link $\LL$, which thus provides a framed link in $Y'$. For each $I \in \{a,b\}$ we have that $Y'(I) \cong Y(I) \# S^{1} \times S^{2}$. We define a map $X(Y,\mathcal{C}) \ra X(Y',\mathcal{C})$ by taking the sum of the maps
${\bf x} \in \widehat{CF}(Y(I)) \ra {\bf x} \times \Theta^{+}_{\frac{1}{2}} \in \widehat{CF}(Y'(I))$, where we implicitly identify $\widehat{CF}(Y'(I)) \cong \widehat{CF}(Y(I)) \otimes \widehat{CF}(S^{1} \times S^{2}, \mathfrak{s}_{0})$. For each theories in which we are interested here, we will have $\HomT{Y(I) \# S^{1} \times S^{2}} \cong \HomT{Y(I)} \otimes \widehat{HF}(S^{1} \times S^{2})$. For the knot Floer homology, the homology for the $S^{1} \times S^{2}$ summand occurs with filtration $0$. \\
\\
\noindent The differential $D'$ for $X(Y',\mathcal{C})$ can be found by noting that in the genus $1$ portion of the diagram representing $S^{1} \times S^{2}$ we always have a $\mu = k-2$ homotopy class of holomorphic discs joining the respective $\Theta^{+}$'s for the different codes. This follows from one of our local models. We may apply the gluing lemma to reduce to the differential from $X(Y,\mathcal{C})$. We will then have $D'(\xi \otimes \Theta^{+}) = (D\xi) \otimes \Theta^{+}$. Thus this map is a (filtered) chain map which induces a map on each spectral sequence. We now deduce several
consequences. The first is straightforward,

\begin{lemma}
The map induced on $E^{1}(X(Y, \mathcal{C}))$ is $\bigoplus_{I \in \mathcal{C}} F_{I}$ where $F_{I}$ is the one-handle map 
$\HomT{Y(I)} \ra \HomT{Y(I) \# S^{1} \times S^{2}}$.
\end{lemma}

\begin{lemma}
The induced one handle map $$H_{\ast}(X(Y,\mathcal{C})) \ra H_{\ast}(X(Y', \mathcal{C}))$$ represents the map $$\widehat{HF}(Y(\{c\ldots c\})) \ra \widehat{HF}(Y'(\{c\ldots c\}))$$ under the quasi-isomorphisms defining the link surgery spectral sequence. 
\end{lemma}

\noindent {\bf Proof:} By the argument above, we actually obtain chain maps for each $X(\mathcal{S})$, including for $\mathcal{S} = \{a,b\}^{n-1} \times \{a,b,c\}$. It is easily verified that these provide a map of short exact sequences:
 
\noindent Altogether, we obtain the following commutative diagram: 
$$
\begindc{\commdiag}[6]
\obj(0,10)[T0]{$0$} 
\obj(15,10)[T1]{$X(Y, \{a,b\}^{n} \times \{c\})$}  
\obj(35,10)[T2]{$X(Y, \mathcal{S})$}  
\obj(55,10)[T3]{$X(Y, \mathcal{C})$} 
\obj(65,10)[T4]{$0$}
\mor{T0}{T1}{$\,$}
\mor{T1}{T2}{$i$}
\mor{T2}{T3}{$\pi$}
\mor{T3}{T4}{$\,$}
\obj(0,0)[B0]{$0$} 
\obj(15,0)[B1]{$X(Y', \{a,b\}^{n} \times \{c\})$}  
\obj(35,0)[B2]{$X(Y', \mathcal{S})$}  
\obj(55,0)[B3]{$X(Y', \mathcal{C})$} 
\obj(65,0)[B4]{$0$}
\mor{B0}{B1}{$\,$}
\mor{B1}{B2}{$i$}
\mor{B2}{B3}{$\pi$}
\mor{B3}{B4}{$\,$}
\mor{T1}{B1}{$F_{0}$}
\mor{T2}{B2}{$F$}
\mor{T3}{B3}{$F_{1}$}
\enddc
$$
\ \\
\noindent In homology this corresponds to having a map of long exact sequences. Since $H_{\ast}(X(Y,\mathcal{S}))$ and $H_{\ast}(X(Y',\mathcal{S}))$ are both trivial, we obtain the following square
$$
\begindc{\commdiag}[6]
\obj(0,10)[I1]{$H_{\ast}(Y, \{a,b\}^{n} \times \{c\})$}
\obj(25,10)[I2]{$H_{\ast}(X(Y, \mathcal{C}))$}
\obj(0,0)[I3]{$H_{\ast}(Y', \{a,b\}^{n} \times \{c\})$}
\obj(25,0)[I4]{$H_{\ast}(X(Y', \mathcal{C}))$}
\mor{I1}{I2}{$\cong$}
\mor{I3}{I4}{$\cong$}
\mor{I1}{I3}{$F_{0\ast}$}
\mor{I2}{I4}{$F_{1\ast}$}
\enddc
$$
The isomorphisms in this square are those induced by the quasi-isomorphisms of chain complexes in the proof of the surgery spectral 
sequence. We can therefore repeat the induction to see that the map induced on the homology of $X(\mathcal{C})$ by the one handle addition is the map
induced on the homology of $Y(\{cc\ldots c\})$, modulo the quasi-isomorphisms. $\Diamond$

\subsection{Maps from three handle attachments}
\ \\
\noindent Suppose now that $W$ is obtained from a three handle attachment to $Y \times I$. As long as $\LL$ misses the attaching sphere in $Y$, this same sphere exists in each of the $Y(I)$. We can thus obtain a Heegaard diagram for $\LL$ which has a genus $1$ summand representing this sphere. This summand exists for all the codes, and so we may repeat the argument given for $1$-handles. In this case, we define a map 
$\widehat{CF}{Y(I)} \cong \widehat{CF}(Y'(I))$ where $\widehat{CF}(Y(I)) \cong \widehat{CF}(Y'(I)) \otimes \widehat{CF}(S^{1} \times S^{2})$. This map takes ${\bf x} \otimes \Theta^{-}_{-\frac{1}{2}} \lra {\bf x}$ and ${\bf x} \otimes \Theta^{+}_{\frac{1}{2}} \lra 0$. That this is a chain map follows from the argument above, as do the following conclusions:

\begin{lemma}
The map induced on $E^{1}(X(Y, \mathcal{C}))$ is $\bigoplus_{I \in \mathcal{C}} F_{I}$ where $F_{I}$ is the three-handle map 
$\HomT{Y'(I)} \otimes \HomT{S^{1} \times S^{2}} \ra \HomT{Y(I)}$.
\end{lemma}

\begin{lemma}
The induced three handle map $$H_{\ast}(X(Y,\mathcal{C})) \ra H_{\ast}(X(Y', \mathcal{C}))$$ represents the map $$\widehat{HF}(Y(\{c\ldots c\})) \ra \widehat{HF}(Y'(\{c\ldots c\}))$$ under the quasi-isomorphisms defining the link surgery spectral sequence. 
\end{lemma}

\subsection{Two handles}
We consider a framed link in $Y$ consisting of two parts, $\LL$ as above, and $\LL'$ which will encode the two handles used in constructing the cobordism. We will denote the result of surgery on $\LL'$ using subscripts; for example, by $Y_{\LL'}$. The components of $\LL$ will be provided
with two framings in a triad, namely the $a$ and $b$ framings. Since we are now dealing with $X(\mathcal{C}, \LL)$, we will change the code to 
be in $\{0,1\}^{n}$. Note that $\LL'$ provides a framed link in $Y(I)$ for any code. When dealing with the twisted coefficient theory or with knot Floer homology we assume that no component in $\LL'$ intersects the spanning surface or the surface determining the twisted coefficients algebraically a non-zero number of times. To
define the chain map on $X(Y, \mathcal{C}, \LL)$ we add to the code above, using $\{0,1\}^{n+1}$, where the first $n$ elements still determine the
framing in the triad. The final element in the code will be interpretted using: $0 \Rightarrow$ use (small Hamiltonian isotopes) of the meridians
of the components in $\LL'$, $1 \Rightarrow$ use (small Hamiltonian isotopes of) the framing curves in $\LL'$ instead of the meridian. Thus 
$Y(I \cup\{1\})$ is another name for $Y_{\LL'}(I)$. Our chain map is then 
$$
F(\xi) = \sum_{J \in \{0,1\}^{n}\times \{1\}}\  \sum_{I = \codeS{0}{k} = J} D_{\codeS{0}{k}}(\xi)
$$
where the second sum occurs over immediate successor sequences starting at $I \in \{0,1\}^{n} \times \{0\}$. Thus as we travel around each polygon, at some point, all the meridians on $\LL'$ flip to their framings {\em simultaneously}, and then we proceed with sequences for $Y_{\LL'}$. Of course, for the different theories we will need to use different maps in place of $D$. These should be adjusted as previously. 

\begin{prop} The map $F: X(Y, \mathcal{C}) \ra X(Y_{\LL'}, \mathcal{C})$ is a filtered chain map.
\end{prop}

\noindent {\bf Proof:} As usual we consider the boundaries of the compactifications of the $1$-dimensional moduli spaces. This provides the $A_{\infty}$-relation noted above. Those boundaries that come from
dividing along the $\al$-curves include the change of the last element in the code from $0$ to $1$ on one side or the other of the divide. The other side corresponds to an immediate successor sequence involving only those codes for $\LL$. Adding over all moduli spaces, these divisions yield the
map $D_{Y_{\LL'}} \circ F + F \circ D_{Y}$. For $F$ to be a chain map, we must see that all of the moduli space boundaries coming from other divisions contribute nothing. For this we need an analog of the cancellation lemma. \\
\ \\
\noindent In particular, we can prove the cancellation as above. For $k > 2$, a $k+1$-gon consisting of $\eta(I)$ labellings on its boundary
either contains a change in the last entry in the code or it doesn't. If it doesn't, then the previous cancellation lemma applies, as the immediate
successor sequence arises either in the spectral sequence for $(Y, \LL)$ or $(Y_{\LL'}, \LL)$. In fact, this applies independently of $k$. However,
when there is a change in code, we must verify the property. The key is to note that the change still occurs in a genus $1$ component, and so we can still use the local model to construct a $k+1$-gon with Maslov index $k-2$. As we only change the framing once in each genus $1$-component, we need
only verify this for the $\hat{\ }$-theory since any other homotopy class with $n_{w} \geq 0$ will be this one plus doubly periodic domains, plus copies of $[\Si]$, and the latter {\em increases} the Maslov index. This leaves us with the $k \leq 2$ case. As before, verifying $k=1$ corresponds to
seeing that the $\Theta$-generators are closed. This leaves $k=2$. \\
\ \\
\noindent There are three different possibilities. The first edge, travelling clockwise, which includes $1$ in the last entry, occurs as the first edge in the triangle for $k=2$. In this case, all three edges occur in spectral sequence for $Y_{\LL'}$. In particular, the code sequence
looks something like $001, 101, 111$. There is another sequence $001, 011, 111$ whose map has the same image, and thus cancels it. The special edge 
can also occur as either the second or the third edge in the triangle. In these cases, we have codes $000$, $100$, $101$ or $000, 010, 011$ for the
last position, and $000, 001, 101$ or $000, 001, 011$ if it occurs in the second. We note that those codes for the last position cancel against those
for the middle position, but not against each other. Nevertheless a local model shows that the two pairs have the same images. As a result, the $k=2$
case also cancels, implying that the map $F$ is, indeed, a chain map. $\Diamond$\\
\ \\
\noindent The top level of this map consists of maps between $I \times \{0\}$ and $I \times \{1\}$. These maps count triangles only, and are easily
seen to be those from \cite{Smoo}. In other words,

\begin{lemma}
The map induced by $F$ on $E^{1}(X(Y, \mathcal{C}))$ is $\bigoplus_{I \in \mathcal{C}} F_{I, \LL'}$ where $F_{I, \LL'}$ is the two-handle cobordism map $\HomT{Y(I \times \{0\})} \ra \HomT{Y(I \times \{1\})}$ found by summing over all $Spin^{c}$ structures on $W$.
\end{lemma}

\noindent There are two examples worth noting. The first is $0$ surgery on an unknot, i.e. a circle bounding a disc, unlinked from the rest of the diagram. This gives a map on $E^{1}$ of $\HomT{Y(I)} \ra \HomT{Y(I)} \otimes \widehat{HF}(S^{1} \times S^{2})$ given by $\xi \ra \xi \otimes \Theta^{-}_{-\frac{1}{2}}$. Unlike for one handle additions, there may be ``higher'' terms in the chain map, i.e. portions of the map with image in lower filtration levels for the flattened filtration. Nevertheless, this $E^{1}$ map is a chain map for the differential on the $E^{1}$-page, and determines the morphism of spectral sequences as usual. In addition, we can consider $-1$ surgery on the same unknot, and the $E^{1}$ map will now be the $0$-map, by the blow-up formula since we sum over all $Spin^{c}$ structures. As a result, the morphism induced on the spectral sequence is also $0$ at every page.  \\
\ \\
\noindent We now turn to showing that the chain map $F$ corresponds to the two handle map $F: \HomT{Y(\{c\ldots c\})} \ra \HomT{Y_{\LL'}(\{c\ldots c\})}$. We do this only for the case when $\mathbb{L}'$ is a knot, and leave the general case for a moment. Let $\gamma \subset Y(\{cc\ldots c\})$ be a framed knot, and let $\mathcal{S}=\{a,b\}^{n-1}\times\{a,b,c\}$. To each code in $\mathcal{S}$ we associate
$Y(I)$, along with the framed knot $\gamma$. We can define a chain map on $X(\mathcal{S})$ by 
$$
F(\xi) = \sum_{J \in \{0\}\times\{a,b\}^{l}\times \{a,b,c\}}\  \sum_{I = \codeS{0}{k} = J} D_{\codeS{0}{k}}(\xi)
$$
where the first entry, $\{0\}$, in the code indicates that we use the framing on $\gamma$. The sequence of immediate successors treats $\gamma$ as another component of $\mathbb{L}$, with the framing defining the triad on $\gamma$. If we let $\overline{\mathbb{L}}$ be the result of adding $\gamma$ to $\mathbb{L}$, then we have a new complex $X(\mathcal{S}')$ for $\mathcal{S} = \{\infty, 0 \} \times \{a,b\}^{n-1}\times \{a,b,c\}$. However, we know that $D^{2} \equiv 0$ on $X(\mathcal{S}')$, and the chain map above is the portion of the differential which changes the code on $\gamma$.\\
\ \\
\noindent Restricting to those codes with a $c$ in the final position, we obtain for $\xi \in \mathcal{S}$, that
$$
F_{0}(\xi) = \sum_{J \in \{0\}\times\{a,b\}^{n-1}\times \{c\}}\  \sum_{I = \codeS{0}{k} = J} D_{\codeS{0}{k}}(\xi)
$$
the chain map for surgery on $\gamma$ for the bouquet derived from $\mathbb{L}$ by $c$ surgery on the final component. Since having $c$ as a final component gives a subcomplex of $X(\mathcal{S})$, with or without surgery on $\gamma$, we can take quotient complexes and look at the resulting map.
The quotient complex for $X(\mathcal{S})$ will be $X(\mathcal{C})$, and the induced map will be
$$
F_{1}(\xi) = \sum_{J \in \{0\}\times\{a,b\}^{n}}\  \sum_{I = \codeS{0}{k} = J} D_{\codeS{0}{k}}(\xi)
$$
which is the chain map on $X(\mathcal{C})$. \\
\ \\
\noindent Altogether, we obtain the following commutative diagram: 
$$
\begindc{\commdiag}[6]
\obj(0,10)[T0]{$0$} 
\obj(15,10)[T1]{$X(\{\infty\}\times\{a,b\}^{n-1} \times \{c\})$}  
\obj(35,10)[T2]{$X(\{\infty\}\times\mathcal{S})$}  
\obj(55,10)[T3]{$X(\{\infty\} \times \mathcal{C})$} 
\obj(65,10)[T4]{$0$}
\mor{T0}{T1}{$\,$}
\mor{T1}{T2}{$f$}
\mor{T2}{T3}{$g$}
\mor{T3}{T4}{$\,$}
\obj(0,0)[B0]{$0$} 
\obj(15,0)[B1]{$X(\{0\}\times\{a,b\}^{n-1} \times \{c\})$}  
\obj(35,0)[B2]{$X(\{0\}\times\mathcal{S})$}  
\obj(55,0)[B3]{$X(\{0\} \times \mathcal{C})$} 
\obj(65,0)[B4]{$0$}
\mor{B0}{B1}{$\,$}
\mor{B1}{B2}{$f$}
\mor{B2}{B3}{$g$}
\mor{B3}{B4}{$\,$}
\mor{T1}{B1}{$F_{0}$}
\mor{T2}{B2}{$F$}
\mor{T3}{B3}{$F_{1}$}
\enddc
$$
\ \\
\noindent In homology this corresponds to having a map of long exact sequences. Since $H_{\ast}(\{\infty\}\times X(\mathcal{S}))$ and $H_{\ast}(X(\{0\}\times\mathcal{S}))$ are both trivial, we obtain the following square
$$
\begindc{\commdiag}[6]
\obj(0,10)[I1]{$X(\{\infty\}\times\{a,b\}^{n-1} \times \{c\})$}
\obj(25,10)[I2]{$H_{\ast}(X(\mathcal{C}))$}
\obj(0,0)[I3]{$X(\{0\}\times\{a,b\}^{n-1} \times \{c\})$}
\obj(25,0)[I4]{$H_{\ast}(X(\{0\}\times\mathcal{C}))$}
\mor{I1}{I2}{$\cong$}
\mor{I3}{I4}{$\cong$}
\mor{I1}{I3}{$F_{\ast}$}
\mor{I2}{I4}{$F_{\ast}$}
\enddc
$$
The isomorphisms in this square are those induced by the quasi-isomorphisms of chain complexes in the proof of the surgery spectral 
sequence. We can therefore repeat the induction to see that the map induced on the homology of $X(Y, \mathcal{C})$ by $F$ can be related
to a map induced on the homology of $Y(\{cc\ldots c\})$. This latter map can be explicitly computed from the description given above for $l=0$. There
we see that it is the standard two handle map for surgery on $\gamma$. \\

\begin{lemma}
The induced two handle map $$H_{\ast}(X(Y,\mathcal{C})) \ra H_{\ast}(X(Y', \mathcal{C}))$$ represents the map $$\widehat{HF}(Y(\{c\ldots c\})) \ra \widehat{HF}(Y'(\{c\ldots c\}))$$ under the quasi-isomorphisms defining the link surgery spectral sequence. 
\end{lemma}

\subsection{Composition and Associativity}

\noindent Of course, we now wish to understand reorderings of handles and how that affects the maps above. The only serious question is in reordering the addition of the two handles. To that end, we consider the composition of the maps for adding along $\LL'$ followed by adding $\LL''$. 

\begin{lemma} \label{lem:assoc}
At the level of the $X(\mathcal{C})$ chain complexes, $$F_{\LL''} \circ F_{\LL'} + F_{\LL' \cup \LL''} = H_{\LL' \cup \LL''} \circ D_{Y} + D_{Y_{\LL' \cup \LL''}} \circ H_{\LL' \cup \LL''}$$
\end{lemma}
\ \\
\noindent As an immediate consequence, we see that the induced morphisms on spectral sequences satisfy $F_{\LL''} \circ F_{\LL'} = F_{\LL' \cup \LL''} = F_{\LL'} \circ F_{\LL''}$ starting at the $E^{1}$-page since the induced maps on this page are sums of Heegaard-Floer cobordism maps, which have this property. 

\begin{cor}
The map $F_{\LL'} : H_{\ast}(X(Y,\mathcal{C})) \ra H_{\ast}(X(Y_{\LL'}, \mathcal{C}))$ corresponds to the cobordism map $F_{\LL'}$ on $\HomT{Y(\{c\ldots c\})}$ under the link surgery quasi-isomorphisms.
\end{cor}

\noindent {\bf Proof:} We know this to be true when $\LL'$ is a knot. By the argument above $F_{\LL''} \circ F_{\LL'} = F_{\LL' \cup \LL''}$ as maps
on homology. This also applies to the induced maps on $\HomT{Y(\{c\ldots c\})}$. Hence the result follows by decomposing $F_{\LL'}$ into a composition of maps for framed knots. $\Diamond$.\\
\ \\
\noindent {\bf Proof of Lemma \ref{lem:assoc}:} To do this
we consider successor sequences with two extra entries in the code, one for the switch along $\LL'$ and the other for the switch along $\LL''$. Thus
our sequences must end $00$ to begin with, switch to $10$ at some point since we attach to $\LL'$ first, and switch again to $11$ when we attach to $\LL''$. The code can change in the first $n$ entries anywhere between these steps, but after the last step must end with $11$ from then on. \\
\ \\
\noindent If we look at the composition of the chain maps we see the emergence of $1$-dimensional moduli spaces of $k+1$-gons. We proceed to analyze these in the now familiar pattern. To be specific, we consider polygons which, as we traverse the boundary clockwise, consist of immediate successor
sequences for the enhanced code, with some edge where we change from $00$ to $10$ at the end, and some other edge where we change from $10$ to $11$
{\em following} at some remove, before returning to the $\al$-edge. The boundary of the compactification of such a moduli space consists of the various divisions possible in the $k+1$-gon. As before, any division involving an $\al$ edge contributes to a composition we understand. Namely, those
divisions which have both special edges on one side or the other contribute $H_{\LL' \cup \LL''} \circ D_{Y} + D_{Y_{\LL' \cup \LL''}} \circ H_{\LL' \cup \LL''}$ where the $H$ map is that counting $0$-dimensional moduli spaces of immediate successor sequences involving a $00$ to $10$ to $11$ change
in the last coordinates, at some positions. Those divisions along $\al$ which divide the two special edges into different polygons contirbute the
composition $F_{\LL''} \circ F_{\LL'}$. We are once again left with understanding divisions in the $\eta$-edges. We have seen that those
divisions involving an immediate successor sequence in $Y$, $Y_{\LL'}$ or $Y_{\LL' \cup \LL''}$ all vanish by the original cancellation lemma. Furthermore, by the variation of this lemma in the section on two handle additions, the contribution of maps where only one edge is divided from the $\al$-edge is also $0$. This leaves one case to consider, when a division occurs that separates portions of both special edges from the $\al$ edge. Once again if $k > 2$ we can find a homotopy clas with $\mu = k-2$ since we can consider the genus $1$ summands separately. We are left with considering the case $k=2$ as $k=1$ is handled as before. For $k=2$, the triangle has three possibilities. Namely, there is some other framing curve
and the triple of codes gives $000$ to $010$ to $011$, or $010$ to $110$ to $111$, or $010$ to $011$ to $111$. In the last two cases, the first time
we have a $1$ in the penultimate spot is the edge in the triangle, or else the cancellation follows from a previous remark. A local picture shows that
the last two sequences will cancel, and the first does not. When the first division occurs, the polygon including $\al$ can thus contribute to a map in which there is an vertex where the edge on one side has $00$ at the end of its code, and the following edge has $11$. In other words, a place where the framings on $\LL' \cup \LL''$ change {\em simultaneously}. If we add over all successor sequences we obtain the associativity result. $\Diamond$

\subsection{Handle Cancellation, Interchange, and Handleslides}

\noindent We now wish to see how the maps above act in relation to each other. Since we are interested principally in the morphisms induced on the spectral sequence we will take the easy route and forgo the chain map calculations. In the preceding subsections we have taken pains to highlight that the maps on the $E^{1}$ page of the spectral sequence are induced by the Heegaard-Floer cobordism maps on the individual pieces $\HomT{Y(I)}$ for each homology theory. We can therefore conclude

\begin{thm}
The cobordism maps on $X(Y,\mathcal{C})$ defined above give rise to natural morphisms of the associated spectral sequences. This
naturality is encoded in saying that the map for $E^{1}$ commutes with the Heegaaed equivalences used in showing independence
of the spectral sequence from the bouquet. 
\end{thm}

\noindent {\bf Proof:} The maps showing independence from the bouquet are filtered with $E^{1}$-terms given by sums of the corresponding maps on the Heegaard Floer homologies of $Y(I)$. In \cite{Smoo}, propositions 4.6, 4.10, and 4.12 show that the cobordism maps, for each index, commute with the Heegaard equivalence maps induced by isotopies, handleslides and stabilizations. In addition, they show that the two handle maps on homology are independent of the bouquet. Thus, the map on the $E^{1}$ pages of the spectral sequence will also have this property. This will be inherited by all the higher pages as well, since whichever composition in the commutation we choose we obtain the same map on the $E^{1}$-page, and it is this map which will induce the maps on the other pages. Furthermore, spectral sequences and their morphisms form a category, so the induced maps are compositions of induced maps on all the pages. $\Diamond$\\
\ \\
\noindent To a cobordism $W$ between $(Y,\mathbb{L})$ and $(Y',\mathbb{L})$ where the handle additions respect $\mathbb{L}$, we can associate
a cobordism map $F_{W}: X(Y,\mathcal{C}) \ra X(Y',\mathcal{C})$ by the composition $F_{3} \circ F_{\LL'} \circ F_{1}$, where $F_{1}$ includes all the
one handles, and $F_{3}$ includes all the three handles. This map is filtered, and we have

\begin{thm}
The map induced on the $E^{1}$ pages of the spectral sequences by $F_{W}$ depends on $W$ and not the particular handle structure chosen
to define $F_{W}$. The same holds for the map induced on each higher page. Under the quasi-isomorphisms between $H_{\ast}(X(Y,\mathcal{C})$ and
$\HomT{Y(\{c\ldots c\})}$, the map $F_{W}$ induces the Heegaard-Floer map $\widehat{F}_{W}$.
\end{thm}

\noindent{\bf Proof:} In \cite{Smoo}, \POz and \ZSz prove the following facts about the Heegaard-Floer cobordism maps {\em on the homology}:
\ben
\item Lemmas 4.13, 4.15: The order of the attatchment of $1$ and $3$ handles does not alter the map induced by a cobordism consisting only of such handles. Nor do handleslides among such handles. 
\item The order of the attatchment of $1$ handles versus two handles which do not go over the one handles does not alter the
map induced on homology. (Implicit in the defintion of the cobordism map)
\item The order of the attachment of $2$ handles versus three handles which do not go over the $2$-handles does not alter the map induced
on homology by the cobordism. Nor does the interchange of one and three handles. (Implicit in the definition of the cobordism map)
\item Lemma 4.14: The map $F_{\LL'}$ induced on homology by attaching two handles to $Y(I)$ does not change alter handleslides among the framed components
of $\LL'$.
\een
As a result, and as a result of the associativity of the previous subsection, the morphisms induced on the spectral sequence by composing the morphisms induced by attaching handles will have the same properties, since they have these properties on the $E^{1}$ page. Furthemore, there are results in \cite{Smoo} about the cancellation of handles:
\ben
\item Lemma 4.16: The standard Heegaard-Floer map on $Y(I)$ induced by adding a one-handle followed by the cobordims map induced by adding a cancelling two-handle
is the identity map.
\item Lemma 4.17: The composition of the cobordism map found by adding a framed two handle to $Y(I)$ followed by a three handle cancelling this two handle is also the identity map.
\een
These depend upon the naturality argument already established in the previous theorem. Again, since the compositions will therefore induce the identity map in the $E^{1}$ page, they induce the identity on all pages of the spectral sequence. $\Diamond$

\section{Duality}

\subsection{Mirror and cohomology}

\noindent We restrict to the hat theories, for now. We can construct a spectral sequence converging to $\mir{Y}$ as well. The link $\mathbb{L}$ defines a link in $\mir{Y}$. The only alteration is in the triad of framings. When we change the orientation on $Y$, we also change the orientation on $\partial Y$, so our triad now has $+1$ intersection for each pair. If we reverse the order $\ga, \be, \al$ we now have $-1$ intersections $\ga \cap \be$ $=-(\be \cap \ga)$ since the intersection form on a surface is anti-symmetric. This also reverses the order of the maps in the triad. We will keep the framing curve marked $c$ in the same position of the cycle. This interchanges the which other framings are marked $a$ and $b$. We wish to examine these maps for the mirror 
$\mir{Y}$. \\
\ \\
\noindent We observe that if $(\Si, \eta_{0}, \ldots, \eta_{k})$ defines a four manifold, $W$, with $k+1$ boundary components, then we can also 
describe $W$ as $(-\Si, \eta_{0}, \eta_{k}, \ldots, \eta_{1})$. This cobordism is $W$ turned upside down, and we will denote it $M$. This arises from flipping the orientation on both the surface and the underlying $k+1$-gon, and using $\mir{H_{\eta_{i}}} \times \mir{I}$ to fill in the $\eta_{i}$ boundary. In addition, if this original Heegaard $n$-tuple is built from an immediate successor sequence for $\eta_{1}, \ldots, \eta_{k}$, then the
new Heegaard tuple is built from an immediate predescesor sequence for the original triad. Once we reverse the order of all triads, we once again obtain an immediate succesor sequence. \\
\ \\
\noindent The key to analyzing the cohomology is that the moduli space for a homotopy class $\phi \in \pi_{2}({\bf x}, z_{1}, \ldots, z_{k-1}, {\bf y})$ in the diagram for $W$ is isomorphic to that for $\pi_{2}({\bf y}, z_{k-1}, \ldots, z_{1}, {\bf x})$ in $M$. Indeed, The family of $J$'s on $Sym^{g}(\Si)$ gives rise to a family of $J$'s on $Sym^{g}(-\Si)$ by taking $J \ra -J$. Then a map $\phi : D \ra Sym^{g}(\Si)$, where $D$ is the 
 $k+1$-gon (thought of as lying in the upper half-space of $\R^{2}$ with the $\eta_{0}$ side on the $x$-axis), can be taken to a map $\phi: D' \ra Sym^{g}(-\Si)$, where $D' = r(D)$ and $r(x,y) = (x,-y)$. If $\phi$ is pseudo-holomorphic for $J$, then $\phi'$ is pseudo-holomorphic for $-J$. \\
\ \\
\noindent However, the maps go in the opposite direction as those for $W$. In other words, we build a model for the cohomology of $Y$. This proves that the chain complex $X(Y, \mathcal{C})$ has cohomology isomorphic to $X(\mir{Y}, \mathcal{C})$. Furthermore, there is a duality pairing between
the part with code $I$ and that with code $\{11\ldots1\} - I$. Since we are working over a field, the pairing is well-defined on the homology of the first page. Furthermore, the induced differentials should also be those for the cohomology since $\mathrm{Hom}(H_{\ast}, \mathbb{F}) = H^{\ast}$ over a field. This is natural for maps between complexes by the universal coefficient theorem, thus we can propagate down through the spectral sequence to obtain a theorem of the form:

\begin{theorem}
$$E^{k}_{\ast}(\mir{Y}, \mir{\LL}) \cong E^{k,\ast}(Y, \LL) \hspace{.5in} E^{k}_{\ast}(\mir{Y}, K, \mir{\LL}, -k) \cong E^{k,\ast}(Y, K, \LL, k)$$ 
\end{theorem}    

\noindent We can then make use of the results in appendix B to conclude that there are duality pairings:
$$
E^{k}_{\ast}(\mir{Y}, \mir{\LL}) \otimes E^{k}_{\ast}(Y, \LL) \ra \F
$$
induced from the duality pairing of $\widehat{CF}(Y, I)$ for any code $I$ with the corresponding group for $\mir{Y}$ and code $I'$ 
where $I'$ is found from $I$ by interchanging the $a$'s and $b$'s. Thus we will pair intersection points ${\bf x}$ in code $I$ and ${\bf y}$ in code $J$, provided that $I$ can be obtained from $J$ by switching the $a$'s and $b$'s, that ${\bf x} = {\bf y}$ as an intersection point for the three manifold corresponding to these codes.  
\ \\
\noindent Suppose we have a framed link $\LL \subset Y_{1}$, as above. For each cobordism $W: Y_{1} \ra Y_{2}$ found by attaching one, two, and three handles, we have constructed a filtered chain map $F_{W} : X(Y_{1}, \LL) \ra X(Y_{2}, \LL)$. Of course $W$ can also
be thought of as a cobordism $W: \mir{Y_{2}} \ra \mir{Y_{1}}$ and this will induce a map $F_{W}': X(\mir{Y_{2}}, \mir{\LL}) \ra X(\mir{Y_{1}}, \mir{\LL})$. Using the argument above, we see that $F_{W}'$ is equivalent to $F_{W}^{\ast}$ (in the notation of the appendix) under the isomorphisms with the cohomology (see also the duality dection in \cite{Smoo}). We can thus conclude:

$$
\langle F_{W,i}(\xi), \eta \rangle_{E^{i}(Y_{2},\LL) \otimes E^{i}(\mir{Y_{2}}, \mir{\LL})} = \langle \xi, F_{W,i}'(\eta) \rangle_{E^{i}(Y_{1},\LL) \otimes E^{i}(\mir{Y_{1}}, \mir{\LL})}
$$

\appendix

\section{Homological Algebra over $\Z$}

\subsection{Mapping Cones}

\noindent Let $f: A \ra B$ be an anti-chain map between chain complexes $(A, \der_{A})$ and $(B, \der_{B})$, i.e. $\der_{B} \circ f + f \circ \der_{A} = 0$. Such a map induces a map on the homologies of the respective chain complexes. Then we can form a new complex, called the mapping cone of $f$ and denoted by $MC(f)$ whose chain group is $A \oplus B$ and whose differential is 
$$
\left( \begin{array}{cc}
\der_{A} & 0 \\
f & \der_{B} \\
\end{array} \right)
$$
This complex fits into a short exact sequence of complexes:
$$
0 \lra (B, \der_{B}) \lra MC(f) \lra (A, \der_{A}) \lra 0
$$
which induces a long exact sequence of homologies, whose connecting homomorphism is that induced by $f$, $f_{\ast}: \HomT{A} \ra \HomT{B}$. Note that the inclusion and the projection are chain maps. If we have an anti-commuting square of anti-chain maps
$$
\csquare{A}{B}{A'}{B'}{f}{g_{1}}{g_{2}}{f'}
$$
anti-commuting up to chain homotopy, i.e. $f' \circ g_{1} + g_{2} \circ f = \der_{B'}\circ H + H \circ \der_{A}$, we can construct a map $M : MC(f) \lra MC(f')$ which, along with $g_{1}$ and $g_{2}$, induces a natural transformation of the short exact sequences in the category formed from chain complexes considered up to homotopy. If the square anti-commutes exactly then the natural transformation exists in the regular category of chain complexes. \\
\ \\
\noindent A {\em quasi-isomorphism} $f: A \ra B$ is a (anti-)chain map inducing an isomorphism on homology. Two chain complexes, $A$ and $A'$, are {\em quasi-isomorphic} if there is a third complex $B$ and quasi-isomorphisms $\phi_{i}: A_{i} \ra B$. This second step is necessary to form an equivalence relation, as a quasi-isomorphism does not typically have an inverse chain map.  \\
\ \\
\noindent{\bf Note:} Over $\F_{2}$ there is no distinction between chain and anti-chain maps. Thus everything we say here works unfettered for
$\F_{2}$-coefficients.

\subsection{The key lemma of \cite{Doub}}

\noindent The following lemma appears in \cite{Doub}. As can be seen in the main body of this paper, in \cite{Doub} it plays a fundamental role in constructing the spectral sequences we are recounting here. \\

\begin{lemma}\cite{Doub}
Let $\{A_{i}\}_{i=0}^{\infty}$ be a set of chain complexes and let 
$$A_{i} \stackrel{f_{i}}{\lra} A_{i+1}\stackrel{f_{i+1}}{\lra} A_{i+2}\stackrel{f_{i+2}}{\lra} A_{i+3}\stackrel{f_{i+3}}{\lra} A_{i+4}\stackrel{f_{i+4}}{\lra} A_{i+5}$$ 
be a sequence of anti-chain maps satisfying the properties:
\ben
\item For each $i$, $f_{i+1}\circ f_{i}$ is chain homotopically trivial, with chain homotopy $H_{i} : A_{i} \ra A_{i+2}$.
\item For each $i$, the map $f_{i+2}\circ H_{i} + H_{i+1} \circ f_{i} : A_{i} \ra A_{i+3}$ is a quasi-isomorphism. 
\een
then $M(f_{i})$ and $M(f_{i+3})$ are quasi-isomorphic to $A_{i+2}$ and $A_{i+5}$, and $M(f_{i+1})$ and $M(f_{i+4})$ are quasi-isomorphic to
$A_{i}$ and $A_{i+3}$. 
\end{lemma}

\noindent {\bf Proof:} Let $$\psi_{i} = -f_{i+2} \circ H_{i} - H_{i+1} \circ f_{i}\ \ :A_{i} \ra A_{i+3}$$ According to the first hypothesis
this map is a chain map. Indeed, $\psi_{i} \circ \der_{i} = -f_{i+2} \circ H_{i} \circ \der_{i} - H_{i+1} \circ f_{i} \circ \der_{i}$ $= -f_{i+2} \circ \big( f_{i+1} \circ f_{i} - \der_{i+2} \circ H_{i} \big) + H_{i+1} \circ \der_{i+1} \circ f_{i}$ $= - f_{i+2} \circ f_{i+1} \circ f_{i} + f_{i+2} \circ \der_{i+2} \circ H_{i} + \big( f_{i+2} \circ f_{i+1} - \der_{i+3} \circ H_{i+1} \big)\circ f_{i}$ $= f_{i+2} \circ \der_{i+2} \circ H_{i} - \der_{i+3} \circ H_{i+1} \circ f_{i}$ $= \der_{i+3} \circ \psi_{i}$. \\
\ \\
\noindent We then note that the following commutes up to homotopy:
$$
\csquare{A_{i}}{A_{i+1}}{A_{i+3}}{A_{i+4}}{f_{i}}{\psi_{i}}{\psi_{i+1}}{f_{i+3}}
$$
Indeed, $\psi_{i+1} \circ f_{i} = -f_{i+3} \circ H_{i +1} \circ f_{i} - H_{i+2} \circ f_{i+1} \circ f_{i}$ $= -f_{i+3} \circ H_{i+1} \circ f_{i} - H_{i+2} \circ \big( \der_{i+2} \circ H_{i} + H_{i} \circ \der_{i}\big)$. But $H_{i+2} \circ \der_{i+2} \circ H_{i} = f_{i+3} \circ f_{i+2} \circ H_{i} - \der_{i+4} \circ H_{i+2} \circ H_{i}$. So $\psi_{i+1} \circ f_{i} = -f_{i+3} \circ H_{i+1} \circ f_{i} - f_{i+3} \circ f_{i+2} \circ H_{i} + \der_{i+4} \circ H_{i+2} \circ H_{i} - H_{i+2} \circ H_{i} \circ \der_{i}$ $= f_{i+3} \circ \psi_{i} + \der_{i+4} \circ H_{i+2} \circ H_{i} - H_{i+2} \circ H_{i} \circ \der_{i}$. So $\psi_{i+1} \circ f_{i} \sim f_{i+3} \circ \psi_{i}$ with anti-chain homotopy $H_{i+2} \circ H_{i}$. In particular, the induced maps on homology will commute. \\
\ \\
\noindent Now we relate to mapping cones by defining two chain maps:
$$
\begin{array}{c}
\al_{i} : M(f_{i}) \lra A_{i+2} \hspace{.5in} \al_{i}(a_{i}, a_{i+1}) =  f_{i+1}(a_{i+1}) - H_{i}(a_{i}) \\
\ \\
\be_{i} : A_{i} \lra MC(f_{i+1}) \hspace{.5in} \be_{i}(a_{i}) = (f_{i}(a_{i}), -H_{i}(a_{i}))
\end{array}
$$
We verify that $\be_{i}$ is an anti-chain map; that $\al_{i}$ is an anti-chain map is similar. $\be_{i} \circ \der_{i} (a_{i}) = (f_{i}(\der_{i}a_{i}), -H_{i}(\der_{i}a_{i}))$ $= (-\der_{i+1}(f_{i}(a_{i})), -f_{i+1}(f_{i}(a_{i})) + \der_{i+2}(H_{i}(a_{i})))$ $= -\der_{MC(f_{i+1})}(f_{i}(a_{i}),$ $-H_{i}(a_{i}))$ $= -\der_{MC(f_{i+1})}\circ \be_{i}(a_{i})$. These maps have the property that $\al_{i+1} \circ \be_{i} = \psi_{i}$ since applied to $a_{i}$, both equal $-H_{i+1}\circ f_{i} (a_{i}) - f_{i+2} \circ H_{i}(a_{i})$. By hypothesis (2), $\psi_{i}$ is a quasi-isomorphism, and thus the composition $\al_{i+1} \circ \be_{i}$ is as well. \\
\ \\
\noindent  Furthermore, $\psi_{i} \circ \pi_{i} : M(f_{i}) \ra A_{i+3}$ is chain homotopic to $f_{i+2}\circ\al_{i}$ by a map $X(a_{i}, a_{i+1}) = H_{i+1}(a_{i+1})$. Note that this is different from \cite{Doub}. Nevertheless, we have $f_{i+2}\circ \al_{i} (a_{i}, a_{i+1}) - \der_{i+3} \circ X(a_{i}, a_{i+1}) - X \circ \der_{MC(f_{i})} (a_{i}, a_{i+1})$ $= -f_{i+2} \circ H_{i} (a_{i}) + f_{i+2} \circ f_{i+1}(a_{i+1}) - \der_{i+3} \circ H_{i+1}(a_{i+1}) - H_{i+1} \circ \big( \der_{i+1}(a_{i+1}) + f_{i}(a_{i})\big)$. Running through the definition of $\psi_{i}$ and $H_{i+1}$ we obtain
$\psi_{i}(a_{i})$ which equals $\psi_{i} \circ \pi_{i}$ on $MC(f_{i})$. \\
\ \\
\noindent A similar argument constructs a chain homotopy between $\be_{i+1} \circ f_{i}: A_{i} \ra MC(f_{i+2})$ and $I_{i+3} \circ \psi_{i}$ using the map $Y(a_{i}) = (H_{i}(a_{i}), 0)$ ($I_{i}$ is the sub-complex inclusion of $A_{i}$ in $MC(f_{i-1})$). \\
\ \\
\noindent The purpose of these calculations is to demonstrate that the squares in the following monster diagram commute up to homotopy.
$$
\begindc{\commdiag}[6]
	\obj(0,0)[I1]{$A_{i+3}$}
	\obj(10,0)[I2]{$A_{i+4}$}
	\obj(20,0)[I3]{$M(f_{i+3})$}
	\obj(30,0)[I4]{$A_{i+3}$}
	\obj(40,0)[I5]{$A_{i+4}$}
	\obj(0,10)[I6]{$A_{i}$}
	\obj(10,10)[I7]{$A_{i+1}$}
	\obj(20,10)[I8]{$A_{i+2}$}
	\obj(30,10)[I9]{$A_{i+3}$}
	\obj(40,10)[I10]{$A_{i+4}$}
	\obj(0,20)[I11]{$A_{i}$}
	\obj(10,20)[I12]{$A_{i+1}$}
	\obj(20,20)[I13]{$M(f_{i})$}
	\obj(30,20)[I14]{$A_{i}$}
	\obj(40,20)[I15]{$A_{i+1}$}
	\mor{I1}{I2}{$f_{i+3}$}
	\mor{I2}{I3}{$I_{i+4}$}
	\mor{I3}{I4}{$\pi_{i+3}$}
	\mor{I4}{I5}{$f_{i+3}$} 
	\mor{I6}{I7}{$f_{i}$}
	\mor{I7}{I8}{$f_{i+1}$}
	\mor{I8}{I9}{$f_{i+2}$}
	\mor{I9}{I10}{$f_{i+3}$}
	\mor{I11}{I12}{$f_{i}$}
	\mor{I12}{I13}{$I_{i+1}$}
	\mor{I13}{I14}{$\pi_{i}$}
	\mor{I14}{I15}{$f_{i}$}
 	\mor{I11}{I6}{$Id$}
	\mor{I6}{I1}{$\psi_{i}$}
 	\mor{I12}{I7}{$Id$}
	\mor{I7}{I2}{$\psi_{i+1}$}
	\mor{I13}{I8}{$\al_{i}$}
	\mor{I8}{I3}{$\be_{i+2}$}
	\mor{I14}{I9}{$\psi_{i}$}
	\mor{I9}{I4}{$Id$}
	\mor{I15}{I10}{$\psi_{i+1}$}
	\mor{I10}{I5}{$Id$}
\enddc 
$$
\ \\
\noindent We have shown that each of the vertical maps is a chain map. When we take homologies, the top and bottom row become the long exact sequences for their respective mapping cones. Since $\psi_{i}$ is a quasi-isomorphism for each $i$, the five lemma applied to the maps between these long exact sequences found by composing the two vertical maps in each column yields that $\be_{i+2} \circ \al_{i}$ is a quasi-isomorphism between the homologies of the mapping cones. On the other hand, $\al_{i+3} \circ \be_{i+2} = \psi_{i+2}$ is a quasi-isomorphism as well. Hence, $\be_{i+2}$ induces a map on homology that is both injective and surjective, i.e. it is a quasi-isomorphism. Thus, $\al_{i}$ and $\al_{i+3}$ induce isomorphisms on homology as well. These induce the first quasi-isomorphisms mentioned in the statement of the lemma. Furthermore, by incrementing all the $i$'s by $1$ in the monster diagram above, we obtain that $\be_{i+3} \circ \al_{i+1}$ is also a quasi-isomorphism. But $\al_{i+1} \circ \be_{i} = \psi_{i}$ is a quasi-isomorphism. Thus $\al_{i+1}$ is a quasi-isomorphism, as are $\be_{i}$ and $\be_{i+3}$.  $\Diamond$\\
\ \\
\noindent {\bf Note:} The middle sequence of our monster diagram yields a complex when we take homology, according to hypothesis (1). The idea of the proof is to squeeze this complex between two long exact sequences. \\
\ \\
\noindent This lemma has an important implication. Let $MC(f_{i}, f_{i+1}, f_{i+2})$ be the {\em iterated mapping cone}, i.e. the chain complex formed from $A_{i} \oplus A_{i+1} \oplus A_{i+2}$ with the differential
$$
\left(
\begin{array}{ccc}
\partial_{i} & 0 & 0 \\
f_{i} & \partial_{i+1} & 0 \\
-H_{i} & f_{i+1} & \partial_{i+2} \\
\end{array}
\right)
$$
That this is a differential follows from the hypotheses in the lemma. There is then a short exact sequence
$$
0 \lra A_{i+2} \lra MC(f_{i}, f_{i+1}, f_{i+2}) \lra MC(f_{i}) \lra 0
$$
The connecting map $\HomT{MC(f_{i})} \lra \HomT{A_{i+2}}$ in the associated long exact sequence is precisely that induced by $\al_{i}$. Since this induced map is an isomorphism by the lemma, we conclude that $\HomT{MC(f_{i}, f_{i+1}, f_{i+2})} \equiv 0$. It is this conclusion which figures centrally in the construction of the Heegaard-Floer spectral sequence.  

\subsection{Key lemma for filtered complexes}

\noindent Let $(A,\mathcal{A})$ and $(B,\mathcal{B})$ be filtered differential modules. Let $f: A \ra B$ be a filtered chain map. Then the mapping
cone $M(f)$ inherits a filtration by declaring $\mathcal{M}_{i} = \mathcal{A}_{i} \oplus \mathcal{B}_{i}$. That the differential preserves this
filtration follows from $f$ being filtered. When undeclared, a filtration on a mapping cone complex will come from this construction.\\

\begin{lemma}
$E^{1}(\mathcal{M}) \cong \mathrm{MC}(E^{1}(f))$
\end{lemma}
\ \\
\noindent {\bf Proof:} $E^{1}(\mathcal{M})$ is the homology of the chain complex whose chain group is $gr\,\mathcal{M}$, i.e. the direct sum $\bigoplus_{i \in \Z} \left(\mathcal{A}_{i}/\mathcal{A}_{i-1} \oplus \mathcal{B}_{i}/\mathcal{B}_{i-1}\right)$ with differential given by 
$$
\left( \begin{array}{cc}
\partial_{A}^{0} & 0 \\
f^{0} & \partial_{B}^{0} \\
\end{array} \right)
$$
Where $f = f^{0} + f^{1} + \ldots$ and $f^{p}$ is that portion of $f$ inducing the modules mapping $\mathcal{A}_{i}/\mathcal{A}_{i-1} \ra \mathcal{B}_{i-p}/\mathcal{B}_{i-p-1}$. Note that this implies that the $E^{1}(\mathcal{M})$ is a mapping cone of the map between the homology of
$\bigoplus_{i \in \Z} \left(\mathcal{A}_{i}/\mathcal{A}_{i-1}\right)$ and the homology of $\bigoplus_{i \in \Z} \left(\mathcal{B}_{i}/\mathcal{B}_{i-1}\right)$ induced by the map $f^{0}$. This map is $E^{1}(f)$ and concludes the lemma. $\Diamond$\\
\ \\
\noindent A filtered chain map $f$ will be called a {\em r-quasi-isomorphism} if it induces an isomorphism between the $E^{r}$ pages of the spectral sequences for the source and the target. For the morphism of spectral sequences induced by $f$, in which the induced maps intertwine the differentials on each page, this implies that all the higher pages, $E^{s}$, are quasi-isomorphic by the induced map, $E^{s}(f)$. Being a 1-quasi-isomorphism is probably weaker than $f$ being a filtered chain isomorphism, but enough for spectral sequence computations.\\
\ \\
\noindent Let $\left\{(A_{i}, \mathcal{A}_{i})\right\}_{i=0}^{\infty}$ be a set of filtered chain complexes with each filtration $\mathcal{A}_{i}$
ascending:
$$
\mathcal{A}_{i} : \{0\} \subset \cdots \subset A^{j}_{i} \subset A^{j+1}_{i} \subset A^{j+2}_{i} \subset \cdots \subset A_{i}
$$
Let $\left\{ f_{i} : A_{i} \ra A_{i+1}\right\}$ be a set of chain maps satisfying:
\ben
\item $f_{i}$ is a filtered map for each $i$.
\item $f_{i+1} \circ f_{i}$ is filtered chain homotopic to $0$, i.e. there is a filtered map $H_{i}: A_{i} \ra A_{i+2}$ such that 
$f_{i+1} \circ f_{i} = \partial_{i+2} \circ H_{i} + H_{i} \circ \partial_{i}$. 
\item Each map $f_{i + 2} \circ H_{i} + H_{i+1}\circ f_{i} : A_{i} \ra A_{i+3}$ is a 1-quasi-isomorphism. 
\een
\ \\
\noindent In this setting we can adapt the key lemma to apply to filtered complexes

\begin{lemma} 
The mapping cone $\mathrm{MC}(f_{2})$ is 1-quasi-isomorphic to $A_{4}$.
\end{lemma}

\noindent {\bf Proof:} The hypotheses above guarantee that the maps in the proof of lemma 4.4 of \cite{Doub} are filtered maps.  We need only check the filtering condition for maps in and out of the mapping cone, but with the aforementioned convention these are clearly filtered. In particular
the map $\psi_{i} = f_{i + 2} \circ H_{i} + H_{i+1}\circ f_{i}$ is a 1-quasi-isomorphism by assumption, and the same argument as in \cite{Doub}
implies that $\alpha_{2}$ is a quasi-isomorphism which is also filtered. This is not quite enough to conclude, but it does ensure that $\alpha_{i}$ induces maps at each page in the spectral sequence.\\
\ \\
\noindent The module $Gr(A_{i}) \cong \oplus_{j\in \Z} A_{i}^{j}/A_{i}^{j-1}$ inherits 
a differential which maps the $j^{th}$ graded component to itself, whose homology provides $E^{1}$. The maps $f_{i}$ induce chain maps between these complexes for each grading level. Indeed each of the maps $\psi_{i}, H_{i}, f_{i}$, etc., likewise induce such maps. Compositions such as $f_{i+1} \circ H_{i}$ induce maps on the graded components which are the same as the compositions for the maps induced from $f_{i+1}$ and $H_{i}$ separately. Thus for each $j$, we have the situation in the lemma in \cite{Doub} applied solely to the $j^{th}$ graded component. Applying the lemma in each grading guarantees that the map induced in that grading by $\alpha_{2}$ is a quasi-isomorphism, i.e. that the induced map on the $E^{1}$ page is an isomorphism of spectral sequences. Thus, $\alpha_{2}$ induces an isomorphism from the $E^{1}$ page for $A_{4}$ to $\mathrm{MC}(E^{1}(f_{2})) \cong E^{1}(\mathrm{MC}(f_{2}))$, which is the desired result. $\Diamond$\\
\ \\
\noindent As in \cite{Doub}, we can reinterpret this as a result on interated mapping cones. Let $M = \mathrm{MC}(f_{1}, f_{2}, f_{3})$ be the 
chain complex defined above for $A_{1} \oplus A_{2} \oplus A_{3}$, filtered by $A_{1}^{j} \oplus A_{2}^{j} \oplus A_{3}^{j}$. The differential defined
above consists of entries that are filtered maps, and thus respects the filtration just defined. The lemma then implies that the induced spectral 
sequence on the iterated mapping cone collapses at the $E^{1}$ term. This follows according to the following diagram:

$$
\begindc{\commdiag}[10]
\obj(0,0)[A00]{$\ $}
\obj(10,0)[A10]{$0$}
\obj(20,0)[A20]{$0$}
\obj(30,0)[A30]{$0$}
\obj(40,0)[A40]{$\ $}
\obj(0,3)[A01]{$0$}
\obj(10,3)[A11]{$A_{3}^{j}/A_{3}^{j-1}$}
\obj(20,3)[A21]{$M^{j}/M^{j-1} $}
\obj(30,3)[A31]{$\mathrm{MC}^{j}(f_{1})/\mathrm{MC}^{j-1}(f_{1})$}
\obj(40,3)[A41]{$0$}
\obj(0,6)[A02]{$0$}
\obj(10,6)[A12]{$A_{3}^{j}$}
\obj(20,6)[A22]{$A_{1}^{j-1} \oplus A_{2}^{j} \oplus A_{3}^{j} $}
\obj(30,6)[A32]{$A_{1}^{j} \oplus A_{2}^{j} $}
\obj(40,6)[A42]{$0$}
\obj(0,9)[A03]{$0$}
\obj(10,9)[A13]{$A_{3}^{j-1}$}
\obj(20,9)[A23]{$A_{1}^{j-1} \oplus A_{2}^{j-1} \oplus A_{3}^{j-1}$}
\obj(30,9)[A33]{$A_{1}^{j-1} \oplus A_{2}^{j-1}$}
\obj(40,9)[A43]{$0$}
\obj(0,12)[A04]{$\ $}
\obj(10,12)[A14]{$0$}
\obj(20,12)[A24]{$0$}
\obj(30,12)[A34]{$0$}
\obj(40,12)[A44]{$\ $}
\mor{A11}{A10}{$\ $}
\mor{A21}{A20}{$\ $}
\mor{A31}{A30}{$\ $}
\mor{A12}{A11}{$\ $}
\mor{A22}{A21}{$\ $}
\mor{A32}{A31}{$\ $}
\mor{A13}{A12}{$\ $}
\mor{A23}{A22}{$\ $}
\mor{A33}{A32}{$\ $}
\mor{A14}{A13}{$\ $}
\mor{A24}{A23}{$\ $}
\mor{A34}{A33}{$\ $}
\mor{A01}{A11}{$\ $}
\mor{A02}{A12}{$\ $}
\mor{A03}{A13}{$\ $}
\mor{A11}{A21}{$\ $}
\mor{A12}{A22}{$\ $}
\mor{A13}{A23}{$\ $}
\mor{A21}{A31}{$\ $}
\mor{A22}{A32}{$\ $}
\mor{A23}{A33}{$\ $}
\mor{A31}{A41}{$\ $}
\mor{A32}{A42}{$\ $}
\mor{A33}{A43}{$\ $}
\enddc
$$
where the top two rows are exact, and all the columns are exact. The nine lemma now guarantees that the bottom row is exact, and each of the maps
is a chain map. In the long exact sequence from the bottom row, there is one map guaranteed to be an isomorphism by the lemma. Consequently, the 
groups in $E^{1}(M)$ are trivial.

\section{Computational Spectral Sequences}

\subsection{Computational Homology} Let $(\mathcal{C}, \partial)$ be a chain complex over a field, $\F$, freely generated by a set $\{x_{i}\}$. That is $\mathcal{C} \cong \oplus \F\,x_{i}$. We review the following standard reduction of such a complex. Suppose $x_{i}$ and $x_{j}$ are two 
generators, such that $\langle \partial x_{i}, x_{j} \rangle = \lambda \neq 0$. Since this is a complex over a field, we may
consider the genrator $x_{i}'=\frac{1}{\lambda} x_{i}$, for which the analogous pairing will equal $1$. We define a new complex
$\mathcal{C}'$ as $\bigoplus \F\,y_{k}$ for $k \neq i, j$ and let $\pi: \mathcal{C} \ra \mathcal{C}'$ be the map determined by the vector space quotient by $\mathrm{Span}\{x_{i}, \partial x_{i}\}$. Since $\partial x_{i} = x_{j} + \sum_{k \neq j} a_{k}x_{k}$, we can think of the quotient as generated by  $y_{k}$ where $y_{k}$ is the image of $x_{k}$ for $k \neq i, j$. This can be equipped with a differential 
$\partial'\,y_{k} = \pi(\partial x_{k})$. Note that written out in the coordinates defined above, this differential becomes 
$$
\partial'\,y_{k} = \pi'(\partial\,x_{k} - \langle \partial x_{k}, x_{j} \rangle \partial\,x_{i})
$$
where $\pi'$ takes $x_{i}, x_{j} \ra 0$ and relabels the other $x$'s to $y$'s. We define another map $\iota : \mathcal{C}' \ra \mathcal{C}$ by $\iota(y_k) = x_{k} - \langle \partial x_{k}, x_{j} \rangle\,x_{i}$. Note that these maps depend upon the choice of cancelling pair. \\

\begin{lemma}
The maps $\pi: (\mathcal{C},\partial) \ra (\mathcal{C}', \partial')$ and $\iota: (\mathcal{C}', \partial') \ra (\mathcal{C},\partial)$ are chain maps.
\end{lemma}
\ \\
\noindent {\bf Proof:} For $\pi$ this follows immediately from from the definition of $\partial'$. For $\iota$ there are three cases to consider: when $y_{k}$ has grading equal to that of $x_{i}$, grading one higher, or a different grading from these. When $y_{k}$ has another grading, then the pairing in the definition of $\iota$ is necessarily $0$, since $x_{j}$ is in grading one less than $x_{i}$. As a result, $\iota(y_{k})$ and $\iota(\partial y_{k})$ both just inject to $x_{k}$ and $\partial(x_{k})$. For elements with the same grading, $\partial\circ \iota\, y_{k}$ is $\partial x_{k} - \langle \partial x_{k}, x_{j}\rangle\,\partial x_{i}$. This is the map $\iota\partial'\,y_{k}$ since $\partial'\,y_{k}$ has the same pattern, and is simply included, for grading reasons, by $\iota$. For an element with grading one higher than $x_{i}$, $\iota(y_{k}) = x_{k}$, and $\partial \circ \iota (y_{k}) = \partial x_{k}$. However, $\partial' (y_{k}) = \pi'(x_{k} - \langle \partial x_{k}, x_{i}\rangle x_{i})$ where we remove the $x_{i}$ component of $\partial x_{k}$. Then $\iota \circ \partial'(y_{k})$ is $\partial x_{k} - \langle \partial x_{k}, x_{i} \rangle x_{i} - \langle \partial(\partial x_{k} - \langle \partial x_{k}, x_{i} \rangle x_{i}), x_{j} \rangle x_{i}$ $= \partial x_{k} - \langle \partial x_{k}, x_{i} \rangle x_{i} + \langle \partial x_{k}, x_{i} \rangle \langle \partial x_{i}, x_{j} \rangle x_{i}$ = $\partial x_{k}$.$\Diamond$\\

\begin{lemma} The maps $\pi$ and $\iota$ are chain homotopy equivalences with 
$$
\begin{array}{c}
\pi \circ \iota - \mathrm{Id} = 0 \\
\ \\
\iota \circ \pi - \mathrm{Id} = \partial H + H \partial
\end{array}
$$
where $H$ is defined by
$$
H(x) = \left\{ \begin{array}{cc} -x_{i} & x=x_{j} \\ 0 & \mathrm{else} \end{array} \right.
$$
\end{lemma}
\ \\
\noindent {\bf Proof:} The map $\pi \circ \iota - \mathrm{Id}$ is identically the zero map, since $\iota(y_{k})$ does not have any coordinate in the direction of $x_{j}$, for any $k$, and the projection cancels the $x_{i}$ coordinates in the image. For $\iota \circ \pi - \mathrm{Id}$, we consider
elements in the same grading as $x_{i}$, $x_{j}$, and otherwise. For those in different gradings, both $\pi$ and $\iota$ are just the maps $x_{k} \ra y_{k}$ and its inverse. Thus $\iota \circ \pi - \mathrm{Id}$ equals $0$, as does the right side above.  For those elements in the same grading as $x_{i}$, the map has the form $-x_{i}$ on $x_{i}$ and $- \langle \partial x_{k}, x_{j} \rangle\,x_{i}$ on the
other elements. This is equal to $H \circ \partial$ for these elements. For those elements in the same grading as $x_{j}$, the map has the form $x_{j} \ra - \sum a_{k}x_{k} - x_{j} = -\partial x_{i}$ and $x_{k} \ra 0$. This is $\partial \circ H$ for these elements, while $H \circ \partial = 0$.$\Diamond$ \\

\begin{thm}
$\mathcal{C}'$ is chain homotopy equivalent to $\mathcal{C}$. \\
\end{thm}

\noindent Now assume that $\mathcal{C}$ is equipped with a filtration, $\mathcal{F}$. Furthermore assume that $\mathcal{F}$ assigns a filtration index to each generator $x_{i}$ along with the subspace it generates. The filtration can then be described as
$$
0 \subset \ldots F_{i} \subset F_{i+1} \subset F_{i+2} \subset \ldots \subset \mathcal{C}
$$
where $\partial(F_{i}) \subset F_{i}$ and $F_{i}/F_{i-1}$ is isomorphic to $\oplus \F\,x_{j}$ for those $x_{j}$ with filtration index $i$. In addition, we can use the filtration indices to decompose the differential:  
$$
\partial = \partial_{0} + \partial_{1} + \partial_{2} + \ldots
$$
where $\partial_{j}$ reduces the fitration index by exactly $j$, i.e. $\partial_{j}\,x = \sum \langle \partial\,x,x_{k}\rangle\,x_{k}$, where
$x$ is homogeneous of degree $i$ and the $x_{k}$ are those generators of degree $i-j$. We now generalize the above theorem.\\

\begin{lemma}\label{lem:special}
Suppose $x_{i}$ and $x_{j}$ such that 
\ben
\item $\partial x_{i} = x_{j} + \sum a_{l} x_{l}$ where $\mathcal{F}(x_{l}) \leq \mathcal{F}(x_{j})$, and
\item if $\langle \partial x_{k}, x_{j} \rangle \neq 0$ then $\mathcal{F}(x_{k}) \geq \mathcal{F}(x_{i})$
\een  
Then $(\mathcal{C}', \mathcal{F}')$ is a filtered chain complex where $\mathcal{F}'(y_{k}) = \mathcal{F}(x_{k})$, and $\pi$ and $\iota$ are filtered chain maps. 
\end{lemma}
\ \\
\noindent {\bf Proof:} Any element with $\langle \partial x_{k}, x_{j}\rangle \neq 0$ occurs in the same level or higher than $x_{j}$. By the first assumption, the alteration in the differential for $y_{k}$ will be by $\sum a_{l}x_{l}$ with the generators in a lower or equal filtration level than that of $x_{j}$. Thus the differential will continue to preserve the filtration, and $\pi$ will be a filtered chain map, since it is the identity on $x_{k}$ for $k \neq i,j$, kills $x_{i}$ and maps $x_{j}$ to $-\sum a_{l}x_{l}$. On the other hand, $\iota(y_{k}) = x_{k} - \langle \partial x_{k}, x_{j} \rangle x_{i}$ will only be filtered if $x_{i}$ has filtration level equal to or lower than that of $x_{k}$. This is ensured by the second assumption. $\Diamond$\\
\ \\
\noindent The assumptions in the theorem clearly apply if we have a filtered complex and $x_{j}$ is in the same filtration level as $x_{i}$. In this case, the chain homotopy map $H$ is also a filtered map. So the reduced complex is filtered chain homotopy equivalent to $\mathcal{C}$. There is another scenario in which the previous theorem applies as well, and we discuss it along with the situation for $H$:\\

\begin{thm}
Let $\mathcal{C}$ be a filtered, freely generated complex over a field with $\partial = \sum_{i \geq k} \partial_{i}$, thereby decreasing the filtration index by at least $k$. Given $x$ and $y$ generators of this complex such that $\langle \partial_{k} x, y \rangle \neq 0$, then the above cancellation produces a filtered complex $\mathcal{C}'$ that is homotopy equivalent to $\mathcal{C}$. Furthermore, the filtered chain maps, $\pi$ and $\iota$, induce isomorphisms of the associated spectral sequences, $E^{r}(\mathcal{C}) \cong E^{r}(\mathcal{C}')$ for $r \geq k+1$ and $r=\infty$. 
\end{thm} 
\ \\
\noindent {\bf Proof:} That the pairing does not vanish implies that $\mathcal{F}(x) - \mathcal{F}(y) = k$. As this is the smallest difference in filtration mediated by a non-zero term in the differential, the assumptions of the previous lemma must hold. To prove the statement about spectral sequences recall that $H$ is zero except on $y$ where $H(y) = -x$ and thus is a chain homotopy of order $k$, since it can increase the filtration index by at most $k$. By Prop. 3.1 in chapter XV of Cartan and Eilenberg, this property of $H$ implies that $\iota\circ \pi$ induces the identity map on $E^{r}(\mathcal{C})$ for $r \geq k+1$. As this is true of $\pi \circ \iota$ for all $r$, we obtain the result. Furthermore, the isomorphism extends
to the $E^{\infty}$ pages as well. Note under the conditions of the theorem prior to this one, we always have a surjective map from $\pi$ and an injective map from $\iota$.$\Diamond$\\
\ \\  
\noindent In the situation of the previous theorem, we also know that $E^{r}(\mathcal{C}) \cong \mathcal{C}$ for $r \leq k$, and this is also true for $\mathcal{C}'$. This allows us to employ the following process to compute the spectral sequence:
\ben
\item Let $\mathcal{C}$ be a filtered chain complex and set $\mathcal{C}_{0} = \mathcal{C}$
\item In $\mathcal{C}_{i}$ find a pair as above where $\langle \partial_{i} x, y \rangle \neq 0$. Since $\partial_{i-1} \equiv 0$, we can
cancel $x$ and $y$ to obtain a new filtered complex, $\mathcal{C}_{i}^{1}$. 
\item Repeat the cancelling of pairs differing in filtration index by $i$, until there are no longer non-zero terms in the differential $\partial_{i}$ for any generator. This gives a sequence $\mathcal{C}_{i}^{j}$ for $j=1, \ldots, N_{i}$. For a finitely supported complex this process always terminates. 
\item Rename $\mathcal{C}_{i}^{N}$ as $\mathcal{C}_{i+1}$. This complex has the property that $\partial_{j} \equiv 0$ for $j \leq i$. Repeat the process starting at step 2, until $\partial \equiv 0$.
\een
We obtain a sequence of filtered complexes $\{\mathcal{C}_{i}\}$ with attendant differentials obtained as above. These complexes are all homotopy equivalent, with differentials that are increasingly long in the filtration indices. Furthemore, the associated graded group for the chain group of $\mathcal{C}_{i}$ is isomorphic to that of $E^{i}(\mathcal{C}_{0})$. In fact, on $\mathcal{C}_{i}$, the induced map $\partial_{i}$ is a differential in itself ($\partial^{2}_{i}\equiv 0$ since all other
terms in $\partial^{2}$ change the filtration value by more than $2i$). The complex $(\mathcal{C}_{i}, \partial_{i})$ has homology
isomorphic to $\mathcal{C}_{i+1}$ since we obtain $\mathcal{C}_{i+1}$ by the cancellation of terms in $\partial_{i}$. We have shown that \\

\begin{thm}
For a filtered, finitely generated complex $\mathcal{C}$ over a field, there is a sequence of homotopy equivalent complexes $\mathcal{C}_{i}$, obtained by the cancellation algorithm above,
such that $\partial_{\mathcal{C}_{i}}$ decreases the filtration by at least $i$, and $E^{i}(\mathcal{C}) \cong (\mathcal{C}_{i}, \partial_{i})$
where $\partial_{i} : \mathcal{C}_{i} \ra \mathcal{C}_{i}$ is the portion of $\partial_{\mathcal{C}_{i}}$ which reduces filtration by exactly $i$.
This sequence stabilizes as $gr H_{\ast}(\mathcal{C})$ after finitely many steps.
\end{thm}
\ \\
\noindent Suppose we have a filtered chain map $f: \mathcal{C} \ra \mathcal{D}$. We can define filtered chain maps $f_{i} : \mathcal{C}_{i} \ra \mathcal{D}_{i}$ by $f_{i}=\pi^{i}_{D}\circ f \circ \iota^{i}_{C}$, where $\pi^{i}$ and $\iota^{i}$ are the maps induced by the cancellation $\mathcal{C} \ra \mathcal{C}_{i}$. Since each of the maps in the composition is a filtered chain map so is $f_{i} : \mathcal{C}_{i} \ra \mathcal{D}_{i}$. These maps induce the same map as $f$ on the associated spectral sequence, starting at $E^{i}$, since the $\pi$ and $\iota$ maps are isomorphisms of spectral sequences starting at this page. Furthermore, it does not matter how we perform the cancellations. We will get maps that are equivalent up to isomorphism of the spectral sequences. This follows since, if $\tilde{\iota}^{i}_{C}, \tilde{\pi}^{i}_{C}$,etc. arise from a cancellation sequence as above, we will have $\pi^{i}_{C} \circ \tilde{\iota}^{i}_{C}$ inducing an isomorphism of spectral sequences for $E^{i}$ and higher. But then $\tilde{\pi}^{i}_{D}\circ \iota^{i}_{D} \circ \pi^{i}_{D} \circ f \circ \iota^{i}_{C} \circ \pi^{i}_{C} \circ \tilde{\iota}^{i}_{C}$ will be chain homotopic to $\tilde{\pi}^{i}_{D} \circ f \circ \tilde{\iota}^{i}_{C}$ by a homotopy of bounded order such that on $E^{i}$ the two maps will induce the same map. However, on $E^{i}$, $\tilde{\pi}^{i}_{D} \circ \iota^{i}_{D}$ is the inverse of $\pi^{i}_{D} \circ \tilde{\iota}^{i}_{D}$, by a bounded chain homotopy, and hence, starting with the $E^{i}$, we must have the same map, up to isomorphism. \\
\ \\
\noindent If $f - g = \partial_{D}\circ H + H\circ\partial_{C}$ where $H$ has order $k$, we can use $\pi^{i}_{D} \circ H \circ \iota^{i}_{C}$ to obtain an order $k$ chain homotopy between $f_{i}$ and $g_{i}$, since the $\pi$ and $\iota$ maps are all chain maps. Starting at $E^{k+1}$, these maps will induce the same map on each additional page of the spectral sequence. \\
\ \\
\noindent It is worth pointing out that the induced map will be $\mathrm{gr}\,f_{i}$ on $E^{i}$. Thus, if the induced map on the spectral sequence is ever an isomorphism, and the filtration is bounded below, we will in fact have isomorphic chain complexes $\mathcal{C}_{k} \cong \mathcal{D}_{k}$ for $k \geq i$.  \\
\ \\
\noindent We note also, as a byproduct of the isomorphism invariance of the maps, that if $\pi$ reduced $\mathcal{D}_{i}$ to $\mathcal{D}_{j}$
and if $\iota$ reverses this map, up to homotopy, for $\mathcal{C}_{j}$, then $\pi \circ f_{i} \circ \iota = f_{j}$, at least up to isomorphism of $\mathcal{C}_{j}$ and $\mathcal{D}_{j}$.

\subsection{Duality}

\noindent To a filtered chain complex as above, we associate another complex $\mathcal{C}^{\ast}$ generated over $\F$ by the same set $\mathcal{G}$.
The filtration index for this complex is given by $\mathcal{F}^{\ast}(x_{i}) = - \mathcal{F}(x_{i})$, and the differential is defined by taking
$$
\partial^{\ast}(x_{i}) = \sum_{x_{j} \in \mathcal{G}} \langle \partial x_{j}, x_{i} \rangle x_{j}
$$
That this is a differential follows from 
$$
\sum_{x_{j} \in \mathcal{G}}  \langle \partial x_{j}, x_{i} \rangle \langle \partial x_{k}, x_{j} \rangle = \langle \partial^{2} x_{k}, x_{i} \rangle = 0
$$
This complex is clearly filtered by $\mathcal{F}^{\ast}$, and arises by ``reversing the arrows'' in the complex $\mathcal{C}$. \\
\ \\
\noindent There is a natural pairing between these complexes, where we denote the generators of $\mathcal{C}^{\ast}$ by $x_{i}^{\ast}$:
$$
\langle x_{j}, x_{i}^{\ast} \rangle = \delta_{ij}
$$
This pairing respects the filtrations in that it is non-zero unless $\mathcal{F}(x_{j}) + \mathcal{F}^{\ast}(x_{i}^{\ast}) = 0$, and has the property that 
$$
\langle \partial x_{j}, x_{i}^{\ast} \rangle = \langle x_{j}, \partial^{\ast} x_{i}^{\ast} \rangle
$$
as is readily verified. The pairing thus descends to a pairing:
$$
H_{\ast}(\mathcal{C}) \otimes H_{\ast}(\mathcal{C}^{\ast}) \lra \F
$$
Out task is to see how the pairing interacts with the spectral sequence. \\
\ \\
We do this by cancellation. For each complex, we have a computational sequence $\{\mathcal{C}_{i}\}$ and $\{\mathcal{C}^{\ast}_{i}\}$. We wish
to see that we can pair these at each stage. To that end, we examine the effect of cancelling a minimal length portion of the differential, $\partial x_{i} = N\,x_{j} + \mathrm{lower\ order}$ in $\mathcal{C}$. Corresponding to this we will cancel $\partial^{\ast} x_{j} = N\,x_{i}^{\ast} + \mathrm{lower\ order}$ in $\mathcal{C}^{\ast}$. The result is two new complexes defined on the same set of generators. We will show that the differentials still have the property defining $\partial^{\ast}$ above. \\
\ \\
\noindent Suppose in $\mathcal{C}$ we have $\langle \partial x_{k} ,x_{j} \rangle = a_{jk}$ and $\langle \partial x_{i}, x_{l}\rangle= a_{li}$. Then,
after the cancellation, we will have $\langle \partial x_{k}, x_{l} \rangle = -\frac{a_{jk}a_{li}}{N}$. If
we apply the same reasoning to $x_{k}^{\ast}, x_{j}^{\ast}, x_{i}^{\ast},$ and $x_{l}^{\ast}$ we will obtain $\langle \partial^{\ast} x_{l}^{\ast}, x_{k}^{\ast} \rangle = - \frac{a_{il}^{\ast}a_{kj}^{\ast}}{N}$. Due to the symmetry of definitions, we have $a_{il}^{\ast} = a_{li}$ and
$a_{kj}^{\ast} = a_{jk}$. Thus we obtain that $\mathcal{C'}^{\ast}$ is the same as $(\mathcal{C}^{\ast})'$. \\
\ \\
\noindent Examining the filtrations shows that the same differential component will have minimal length in each, so if we cancel the same generators at each step we obtain that $(\mathcal{C}_{i})^{\ast} \cong (\mathcal{C}^{\ast})_{i}$. In particular, there is a pairing at each stage in the sequence, {\em and the pairing is that inherited by the remaining generators}. Furthermore, once we have cancelled all the differentials of length $< l$, $\partial_{l}$ is a differential in itself, and 
$$
\langle \partial_{l} x_{i}, x_{j}^{\ast} \rangle = \langle x_{i}, \partial^{\ast}_{l} x_{j} \rangle
$$
Thus the pairing induces a pairing between the homologies $H_{\ast}(\mathcal{C}_{l}, \partial_{l}) \cong E^{l+1}(\mathcal{C})$ and $H_{\ast}(\mathcal{C}^{\ast}_{l}, \partial^{\ast}_{l}) \cong E^{l+1}(\mathcal{C}^{\ast})$. This pairing is the one inherited by $\mathcal{C}_{l+1}$ and $\mathcal{C}_{l+1}^{\ast}$ under the isomorphism found by applying the reduction process to all differentials of length $l$. Altogether,

\begin{prop}
The pairing of $\mathcal{C} \otimes \mathcal{C}^{\ast} \ra \F$ above induces pairings:
$$
E^{l}(\mathcal{C}) \otimes E^{l}(\mathcal{C}^{\ast}) \lra \F
$$
$$
H_{\ast}(\mathcal{C}) \otimes H_{\ast}(\mathcal{C}^{\ast}) \lra \F
$$
\end{prop}
\ \\
\noindent Now suppose we have a filtered chain map $F: \mathcal{C} \ra \mathcal{D}$. Then there is an induced filtered map $F^{\ast} : \mathcal{D}^{\ast} \ra \mathcal{C}^{\ast}$ found by duality:
$$
\langle F^{\ast}(y_{i}^{\ast}), x_{i}^{\ast} \rangle = \langle y_{i}, F(x_{i}) \rangle
$$
$F^{\ast}$ is a chain map. This is seen most easily by noting that $\mathcal{C}^{\ast} \cong \hom (\mathcal{C}, \F)$, and $F^{\ast}$ is the dual map. 
We can form two new complexes, the mapping cones for these maps, namely $\mathcal{C} \oplus \mathcal{D}$ and $\mathcal{C}^{\ast} \oplus \mathcal{D}^{\ast}$ with differentials:
$$
\left(
	\begin{array}{cc}
		\partial_{\mathcal{C}} & 0 \\
		F & -\partial_{\mathcal{D}} \\
	\end{array} 
\right) \hspace{.5in}
\left(
	\begin{array}{cc}
		\partial_{\mathcal{D}^{\ast}} & 0 \\
		F^{\ast} & -\partial_{\mathcal{C}^{\ast}} \\
	\end{array}
\right)
$$
The key observation is $\mathrm{MC}(F^{\ast}) = \mathrm{MC}(F)^{\ast}$. If we cancel the portions of the differentials of length $\leq l$, using
the same sequence for both complexes, and limiting the cancellations to those internal to $\mathcal{C}$, $\mathcal{D}$ and their duals, we obtain
a new set of chain complexes, $\mathrm{MC}(F_{i})$ where $F_{i} : \mathcal{C}_{i} \ra \mathcal{D}_{i}$, and $\mathrm{MC}(F_{i}^{\ast})$. This in turn
implies that $\mathrm{MC}(F_{i}^{\ast})\cong \mathrm{MC}(F_{i})^{\ast}$ and 
$$
\langle F_{i}(x_{i}), y_{i}^{\ast} \rangle_{\mathcal{D} \otimes \mathcal{D}^{\ast}} = \langle x_{i}, F^{\ast}_{i}(y_{i}^{\ast}) \rangle_{\mathcal{C} \otimes \mathcal{C}^{\ast}}
$$
Furthermore, the portion of $F_{i}$ which does not change the filtration index commutes with $\partial_{i}$ and thus produces a chain map
on $(\mathcal{C}_{i}, \partial_{i})$. Following through the implications implies that the induced maps on the spectral sequence 
$$
\langle F_{i,\ast}(x_{i}), y_{i}^{\ast} \rangle_{E^{i}(\mathcal{D}) \otimes E^{i}(\mathcal{D}^{\ast})} = \langle x_{i}, F^{\ast}_{i,\ast}(y_{i}^{\ast}) \rangle_{E^{i}(\mathcal{C}) \otimes E^{i}(\mathcal{C}^{\ast})}
$$

\end{document}